\documentclass[12pt]{amsart}
\usepackage{amsmath,amssymb,amsfonts,a4wide}
\usepackage{epsfig}

\def\cal{\mathcal}

\numberwithin{equation}{section}

\newtheorem {thm}{Theorem}
\newtheorem {lem}[thm]{Lemma}
\newtheorem {pr}[thm]{Proposition}

\def\lbl{\label}
\def\be{\begin{equation}}
\def\ee{\end{equation}}
\def\lbl{\label}
\def\qed{\sq}

\def\la{\lambda}
\def\Bin{\rm Bin}
\def\Geom{\rm Geom}
\def\qed{\hfill $\Box$ \hfill \\}
\def\pf{\noindent{\em Proof:\ }}
\def\var{{\rm var}}

\def\E{{\mathbb E}}
\def\P{{\mathbb P}}
\def\eps{\epsilon}

\def\gg{{\gamma}}
\def\aa{{\alpha}}
\def\bb{\beta}

\begin{document}

\title[Perpetuities with thin tails]{Perpetuities with thin tails, revisited}
\author[Pawe{\l} Hitczenko]{Pawe{\l} Hitczenko$^\dag{}$}
\address{Pawe{\l} Hitczenko\\
Departments of Mathematics and Computer Science  \\
Drexel University\\
Philadelphia, PA 19104 \\
U.S.A}
\email{phitczenko@math.drexel.edu}
\urladdr{
http://www.math.drexel.edu/$\sim$phitczen}

\thanks{$^\dag{}$ This author is supported in part by the NSA grant  \#H98230-09-1-0062}

\author[Jacek Weso{\l}owski]{Jacek Weso{\l}owski 
}
\address{Jacek Weso{\l}owski\\
Wydzia{\l} Matematyki i Nauk Informacyjnych\\
Politechnika Warszawska\\
Plac Politechniki~1\\ 
00-661 Warszawa, Poland}
\email{wesolo@mini.pw.edu.pl}

\date{\today}
\subjclass{Primary:  60H25; secondary: 60E99 }
\keywords{perpetuity, stochastic difference equation, tail behavior}

\maketitle

\begin{abstract}
We consider the tail behavior of random variables $R$ which  are
solutions of the distributional equation $R\stackrel d=Q+MR$,
where $(Q,M)$ is independent of $R$ and $|M|\le 1$. Goldie and
Gr\"ubel showed that the tails of $R$ are no heavier than
exponential and that if   $Q$ is bounded and $M$ 
resembles near 1  the uniform distribution, then the tails of
$R$ are Poissonian. In this paper we further investigate the
connection between the tails of $R$ and   the behavior of $M$ near
1.  We focus on the special case when $Q$ is constant and $M$ is
non--negative.
\end{abstract}

\section{Introduction}
In this note we consider a random variable $R$ given by the
solution of the stochastic equation \be R\stackrel
d=Q+MR,\lbl{perp}\ee where $(Q,M)$ are independent of $R$ on the
right-hand side. Under  suitable assumptions on $(Q,M)$  one can
think of $R$ as a limit in distribution of the following iterative
scheme \be\lbl{perp-nq} R_n=Q_n+M_nR_{n-1},\quad n\ge1\ee where
$R_0$ is arbitrary and $(Q_n,M_n)$, $n\ge1$, are  i.i.d. copies of
$(Q,M)$, and $(Q_n,M_n)$ is independent of $R_{n-1}$. Writing out
the above recurrence and renumbering the random variables
$(Q_n,M_n)$ we see that $R$ may also be defined
 by
\be\lbl{perp-it}R\stackrel d=\sum_{j=1}^\infty
Q_j\prod_{k=1}^{j-1}M_k,\ee provided that the series above
converges  in distribution. Sufficient conditions for the almost
sure convergence are known and have been given by Kesten
\cite{kesten}  who also considered a multidimensional case when
$M$ is a matrix and $Q$ a vector. For a nice detailed discussion
of a one dimensional case  we refer to the paper by Vervaat
\cite{vervaat}; we only mention briefly here that
$\E\log^+|Q|<\infty$ and $\E\log|M|<0$ suffice for the almost sure
convergence of the series in (\ref{perp-it})

In the form (\ref{perp-it}) $R$ has been studied in insurance
mathematics under the name perpetuity. Since schemes like
(\ref{perp-nq}) are ubiquitous in many areas of applied
mathematics, the properties of $R$ have attracted a considerable
interest. We refer to \cite{cl,eg,goldie,gg,kesten,letac,vervaat}
and references therein for more information and sample of
applications. For examples of more recent work on perpetuities and
their applications see  \cite{air,bw, hm,kn}.  A few additional
situations in which perpetuities arise will be mentioned below.

The main focus of research is the tail behavior of $R$. Kesten
\cite{kesten} showed that if $\P(|M|>1)>0$ then $R$ is always
heavy--tailed. More precisely, he showed that if there exists a
$\kappa$ such that $\E|M|^\kappa\log^+|M|<\infty$,
$\E|Q|^\kappa<\infty$, and $\E|M|^\kappa=1$ then for some constant
$C$
$$\P(|R|\ge x)\sim Cx^{-\kappa},\quad\mbox{as}\quad   x\to\infty.$$
Here, and throughout the paper the symbol $f(x)\sim g(x)$ means
that the ratio goes to 1 as $x\to\infty$. His result was
rediscovered, reproved, and extended  by several authors (see
\cite{goldie,grey,grin}). In the complementary case, $\P(|M|\le 1)$
the picture is much less clear. The main  work  we are aware of is
that of Goldie and Gr\"ubel \cite{gg} who showed that in that
case, the tails are never heavier than exponential and that if $M$
behaves near 1 as a uniform random variable then the tails have
Poissonian decay. In their arguments Goldie and Gr\"ubel relied on
inductive arguments applied to (\ref{perp-nq}).

The main purpose of this note is to use systematically their
approach to obtain additional information on the links between the
behavior of $M$ near 1 and the tail behavior of $R$. Following
Goldie and Gr\"ubel (and also customs in large deviation theory)
we will be interested in the asymptotics of the logarithm of the
tail probability, i.e. $\ln\P(|R|\ge x)$ as $x\to \infty.$ Since
we are mainly interested in establishing the links between $M$ and
$R$, we will often make additional, but common, assumptions when
necessary. For example, we generally assume that $Q$ and $M$ are
independent or even that $Q\equiv q$ is non--random. The
independence assumption is typically  needed only for the lower
bounds  on the log of the tail probability, the upper bounds are
usually obtainable without it.  Once the independence of $Q$ and
$M$ is assumed the restriction that $Q$ is degenerate does not
seem to be a major restriction, but makes some of the arguments
more transparent. It is rather the assumption that $Q$ is bounded,
which seems to play the more important role. Similarly, we will
assume that $M$ and $q$ are non--negative.  How the non--negative
case differs from the general is relatively well understood (see
e.g. arguments in \cite[Theorem~2.1, Theorem~3.1, Lemma~5.3]{gg})
to see how arguments for non--negative case can be extended to
more general situations.

We would like to mention an interesting connection of perpetuities
with  a subclass of  infinitely divisible laws, namely, as was
shown by Jurek  \cite{jurek} all self--decomposable random
variables (we refer to \cite{jurek} for the definition) can be
represented as perpetuities $R$ given by (\ref{perp}) with $0\le
M\le 1$. As a matter of fact, much more is shown in \cite{jurek}, namely, if $R$ is
self--decomposable then for {\em  every} random variable $0\le
M\le 1$  there exists a random variable $Q$ (typically not
bounded) such that (\ref{perp}) holds with $(Q,M)$  independent of
$R$ on the right--hand side. This curious result seems to be of
little help as far as general theory of perpetuities goes. In
fact, one can take $M$ to be  any constant $M=m\in(0,1)$ and
equally well represent a self--decomposable random variable as a
series of weighted i.i.d. random variables, with weights forming a
geometric progression. Nonetheless, we mention that building on an
earlier work of Thorin \cite{thorina, thorinb}, Bondesson
\cite{bond} proved a general result which implies, in particular,
that all gamma, inverse gamma, Pareto, log--normal, and Weilbull
distributions are self--decomposable. Some of these results were
obtained earlier by other authors and we refer to Bondesson
\cite[Section~5]{bond} for credits and more examples.

\section{General outline}
To begin the discussion, assume that $ |M|\le 1$. Trivially, if
$|Q|\le q$ and $|M|$ is concentrated on a proper subinterval
$(0,1-\delta)$, $\delta>0$ of $(0,1)$ then the perpetuity $R$ is a
random variable whose absolute value is bounded by $q/\delta$ and
thus has a trivial tail in the sense that $\P(|R|\ge x)=0$ for
$x>q/\delta$. On the other hand if $M$ is not bounded away from 1
then we have the following observation due to Goldie and Gr\"ubel:
\begin{pr}\lbl{prop:lbdd}
For $\delta\in(0,1)$ let $p_\delta:=\P(1-\delta\le M\le1)$. Then
for every such $\delta$ and for all $y>0$ we have
\be\lbl{lbdd}
\P(R\ge\frac q\delta(1-(1-\delta)^y))\ge
p_\delta^{y}.
\ee
In particular, if for $c\in (0,1)$ and $x>q$ we
set
$$
\delta=\frac{cq}x\quad\mbox{and}\quad y=\frac{\ln(1-c)}{\ln(1-cq/x)},
$$
then we get that \be\lbl{lbdd-x} \P(R\ge
x)\ge\left(p_{\frac{cq}x}\right)^{\frac{\ln(1-c)}{\ln(1-cq/x)}}=
\exp\left(\frac{\ln(1-c)}{\ln(1-cq/x)}\ln (p_{\frac{cq}x})\right).
\ee
\end{pr}
\pf
This was observed by  Goldie--Gr\"ubel: For a given $\delta>0$
we let
$$
\tau=\tau_\delta=\inf\{n\ge1:\ M_n<1-\delta\}.
$$
Then by non--negativity and (\ref{perp-it}), on $\{\tau\ge n\}$ we have
$$
R\ge \sum_{k=1}^{n}q(1-\delta)^{k-1}=\frac q\delta\left(1-(1-\delta)^n\right).
$$
Therefore,  for all $n\ge1$,
$$
\P(R\ge\frac q\delta(1-(1-\delta)^n))\ge\P(M_k\ge1-\delta,\ 1\le k<n)= p_\delta^{n-1}.
$$
Hence,
$$
\P(R\ge\frac q\delta(1-(1-\delta)^y))\ge p_\delta^{y},\quad\mbox{for all $y>0$}
$$
which proves (\ref{lbdd});  (\ref{lbdd-x}) follows by a simple calculation.\qed

It is clear from the above proposition  that if $p_\delta$ is
strictly positive for every $\delta>0$ then the perpetuity $R$ has
non-trivial tails.  It is then the behavior of 
$M$ near 1 that  determines the nature of the tails of $R$. It
appears that essentials of such a behavior are shared by a class
of equivalent distributions in the following sense.

Let $\mu$ and $\nu$ be probability distributions on $[0,1]$. For
any $\delta\in(0,1)$ we denote $\mu_{\delta}=\mu((1-\delta,1])$ and
$\nu_{\delta}=\nu((1-\delta,1])$ We say that the distributions 
$\mu$ and $\nu$ are equivalent at 1 if
$$
\exists\  \varepsilon>0\quad\mbox{and}\quad
0<d<D<\infty\quad\mbox{such that} $$ \begin{equation}\label{roro}
\forall\delta\in(0,\varepsilon]:\quad d\le
\frac{\mu_{\delta}}{\nu_{\delta}}\le D\;.
\end{equation}

As we mentioned earlier, our goal here is to shed some additional
light on the relationship between the behavior of the
distribution of $M$ in the left neighborhood of 1 and the tails of
$R$. To accomplish that we will develop in a systematic way the
approach of Goldie and Gr\"ubel. For the upper bound this approach
relies on iteration of (\ref{perp-nq}) to get a uniform upper
bound on the moment generating function of $R_n$ for all $n\ge1$
and then use exponentiation and Markov inequality to translate
this bound into bounds on the tails.  We will develop this in the
next section, but to give a flavor of this argument we provide the
following illustration: Consider (\ref{perp}) and assume that $Q$,
$M$, and $R$ on the right-hand side of (\ref{perp}) are
independent (that is of course stronger than the usual assumption
that $(Q,M)$ are independent of $R$). Also, assume  that $0\le
M\le1$ and that $m:=\E M<1$. To get an upper bound on the moment
generating function $\E e^{zR}$ of $R$, the principle of what
Goldie-Gr\"ubel did is the following:  for $n\ge1$   we have
$$
\E e^{zR_n}=\E e^{z(Q_n+M_nR_{n-1})}=\E e^{zQ}\E e^{zMR_{n-1}}\le
\E e^{zQ}\left\{1+m \E(e^{zR_{n-1}}-1)\right\},$$ where in the last
step we use the fact that for  $s>0$ \be\lbl{dom}\E e^{sM}\le
\E e^{s\Bin(1,m)}=1+m (e^s-1).\ee To set up an induction we seek a
function $A(z)$ such that
\begin{itemize}
\item[(i)]{} $Ee^{zR_{n-1}}\le A(z)$, and \item[(ii)]{}
$\E e^{zQ}\left\{1+m(A(z)-1)\right\}\le A(z)$.
\end{itemize}
Solving (ii) gives
$$B(z):=\frac{(1-m)\E e^{zQ}}{1-m \E e^{zQ}}\le A(z),$$
for $z$ such that $m \E e^{zQ}<1$. Now, $B(z)$ is recognized as the
moment generating function of $\sum_{k=1}^NQ_k$ where $N\stackrel
d=\Geom(1-m)$ and is independent of the sequence $Q_k$, $k\ge1$.
So if we start with any $R_0$ for which (i) holds with $B(z)$ in
place of $A(z)$ then the induction goes through and, under a
reasonably weak assumptions on $Q$, we get an exponential upper
bound on the tail of $R$. In particular if we take $Q\equiv1$ and
$M\stackrel d=\Bin(1,m)$ then $R$ has moment generating function
bounded by that of a geometric random variable and hence
sub-exponential tails as was already shown by Goldie and Gr\"ubel.

We  mention briefly that  the sums described by $B(z)$ are yet
another example of perpetuities. Sums like these  are of interest in renewal theory and risk assessment, for
example. They have been studied before, for instance in
\cite{brown, yan},  under the name
geometric convolutions and geometric random sums, respectively.
 We refer the interested reader there for more information and further  references.

As for the lower bound, the best that is available at this point
is argument based on Proposition~\ref{prop:lbdd}. Interestingly,
this proposition provides a surprisingly good lower bound. By this
we mean the  fact that if the upper bound obtained by  the above
method is constructed carefully so as to be  relatively tight,
then one can usually obtain a  lower bound of a similar strength
from Proposition~\ref{prop:lbdd}. This will be seen in several
situations below. It is thus important to understand how to
construct a tight upper bound. Although we do not have a general
result to that effect, in the last section we will provide
an argument in a particular example that provides a heuristic 
which should  work well in other cases.

The rest of the paper is organized as follows, in the next section
we will discuss an upper bound and in particular, we will state an
inequality (see (\ref{iteration}) below) that is crucial for the
inductive argument. In subsequent sections we will illustrate this
with several examples. Those include beta$(\alpha,\beta)$
densities, and what (for the lack of a better name) we call the
generalized beta$(1,\beta)$ densities. The reason for considering
beta distributions is that one might reasonably hope that they
provide a natural parametrization of a behavior of $M$ near 1,
which could be translated to the tail behavior of $R$. This,
however, is not the case, since as we will show all beta
distributions lead to the same, namely Poissonian, behavior. It
turns out that a much more rapid than power--type variability of
$M$ at 1  is needed to observe a different tail behavior of $R$.
We will then  construct densities for which the logarithm of the
tail probability will have power behavior  $-x^r$, for
$1<r<\infty$. In the last section we will discuss one more example
mainly to illustrate a techinque of constructing $M$ that would give a particular tail behavior of
$R$ in other situations.

\section{Upper bounds}
We begin with the following well--known fact.
\begin{pr}
Suppose that \be\lbl{upbdd-mgf}\E e^{zX}\le\exp(B\Phi(z)),\ee for
some function $\Phi:\ [0,\infty)\to[0,\infty)$, $B>0$ and all
$z>0$. Then \be\lbl{upbdd-tail}\P(X\ge x)\le e^{-\Phi^*(x)}, \ee
where $\Phi^*=\Phi_B^*$ is  defined by
\be\lbl{comple}\Phi^*(x)=\sup\{zx-B\Phi(z):\  z>0\}.\ee
\end{pr}
Note  that if $\Phi$  is an Orlicz function (a convex, continuous,
non--decreasing  function,  such that  $\Phi(0)=0$ and
$\Phi(t)\to\infty$ as $t\to\infty$) then $\Phi^*$ is just a
complementary function to $\Phi$.

\pf This is well--known;
by the usual exponentiation and Markov's inequality we have
$$\P(X\ge x)=\P(e^{zX}\ge e^{zx})\le e^{-zx}\E e^{zX}\le e^{-zx}e^{B\Phi(z)}=e^{-(zx-B\Phi(z))}.
$$
Since the right--hand side may be minimized over $z$ we obtain (\ref{upbdd-tail})
as required.\qed

One can obtain a bound on the moment generating function of $R$
using the fact that it is a limit in distribution of  the
iterative procedure (\ref{perp-nq}) and verifying
(\ref{upbdd-mgf}) for every $R_n$.  In the case $Q_n\equiv q$
(\ref{perp-nq}) takes the form \be\lbl{perp-n} R_n\stackrel
d=q+M_nR_{n-1},\ee where $M_n$ is a copy of $M$ independent of
$R_{n-1}$. To argue inductively, suppose that for some $B>0$
\be\lbl{mgf-n}\E e^{zR_{n-1}}\le \exp(B\Phi(z)),\quad z>0.\ee Then
by (\ref{perp-n}) and (\ref{mgf-n}) applied conditionally on
$M_{n}$ we have
$$
\E e^{zR_n}=e^{qz}\E e^{zM_nR_{n-1}}\le
e^{qz}\E e^{B\Phi(zM_n)}.
$$
The inductive step will be complete once we show that
$$e^{qz}\E e^{B\Phi(zM)}\le e^{B\Phi(z)}.$$
In terms of the distribution $\mu$ of $M$, the above inequality
reads \be\lbl{iteration}e^{qz}\int_0^1 e^{B\Phi(zt)}\mu(dt)\le
e^{B\Phi(z)}.\ee Once this inequality is established, the
induction is complete as one can start with arbitrary random
variable $R_0$, so in particular we can ensure that (\ref{mgf-n})
holds for $R_0$.  The above inequality is crucial for establishing
the upper bound. 

We will be interested in the tail bounds for
large values of $x$. We assume that $\Phi$ is non-degenerate
($\Phi(t)\ne0$ for $t\ne0$) and satisfies $\Phi(t)/t\to\infty$ as
$t\to\infty$ (i.e. $\Phi$ is an $N$--function in the language of
\cite{kr}). Then  $\Phi^*$ has the same properties and  it follows
directly from the definition (\ref{comple}) that as $x\to\infty$
the supremum in (\ref{comple}) is attained at $z\to\infty$. This
means that it suffices that (\ref{upbdd-mgf}) and thus
(\ref{iteration}) hold only for large values of $z$.  Thus, we
have the following consequence of the above discussion:
\begin{pr}\lbl{prop:upbdd}
Let $R$ be given by (\ref{perp}) with $Q\equiv q$. Suppose that
 there exist $B>0$ and $z_0$ such that  %
(\ref{iteration})  is satisfied for the distribution of  $M$ for all $z\ge z_0$. Then
\be\lbl{upbdd-ls}\limsup_{x\to\infty}\frac{\ln \P(R\ge
x)}{\Phi^*_B(x)}\le -1.\ee
\end{pr}

\section{Beta distributions}

As earlier we will denote by $\mu$ the distribution of $M$.
Goldie--Gr\"ubel \cite[Theorem~3.1]{gg}  showed that if $Q$ is
bounded and $\mu$ and the uniform distribution on $[0,1]$ are
equivalent at 1 in the sense of (\ref{roro}) then the resulting
perpetuity has Poissonian tails, that is
$$
\lim_{x\to\infty}\frac{\ln\P(R\ge x)}{x\ln x}=-\frac1q.
$$
Note that uniform and beta $\beta(\alpha,1)$ distributions are
equivalent at 1. One might reasonably hope that considering other
values of the second parameter of the beta distribution might lead
to a different tail behavior of $R$ but this is not the case. As
we show below {\em any}  $M$ whose distribution is equivalent at  1 to a measure with polynomial density at  1 leads to the Poissonian tails of $R$.

\begin{thm} \label{thm:beta}
Let the distribution of $M$ and the $\mbox{beta}(\aa,\bb)$
distribution be equivalent at 1. Assume that $Q\equiv q>0$. Then
$$\lim_{x\to\infty}\frac{\ln\P(R\ge x)}{x\ln x}=-\frac\bb q.$$
\end{thm}
\pf Note that all beta distributions with the same $\bb$ parameter
and different $\aa$ parameters are equivalent in the sense of
(\ref{roro}).  Consequently,  we assume for convenience that
$\aa=1$ so that we consider the beta distribution with the density
$$f(t)=\bb(1-t)^{\bb-1},\quad 0<t<1\;,$$
which is equivalent to the distribution of $M$ at 1.

We show that regardless of the value of $\bb>0$ the tails of the
resulting perpetuities are Poissonian. To get an upper bound we
verify that (\ref{iteration}) holds with $\Phi(z)=e^{bz}$ for a
suitable constant $b$ and some $B>0$. Once this is done, it
follows from the discussion in the previous section that
$$
\ln\P(R\ge x)\le -\frac x{b}\ln\left(\frac x{Bbe}\right)=-\frac1b
x\left(\ln x-\ln(Bbe)\right).$$ which implies that
\be\lbl{betatail}\limsup_{x\to\infty}\frac{\ln\P(R\ge x)}{x\ln
x}\le-\frac1b.\ee Thus we are to show that for sufficiently large
$z>0$ \be\lbl{uppbeta}
e^{qz}\int_0^1\exp(Be^{bzt})\mu(dt)\le\exp(Be^{bz}),\ee for some
positive constant $B$ and $b=q/\bb$. To that end take an $\varepsilon$ for which 
(\ref{roro}) holds with $\nu$ being a $\mbox{beta}(1,\bb)$ distribution. Assume  a  $t_0$ is chosen so that $t_0>1-\varepsilon$. We split the integral on the 
left--hand side as
$$
 e^{qz}\int_0^{t_0}\exp(Be^{bzt})\mu(dt) +
e^{qz}\int_{t_0}^1\exp(Be^{bzt})\mu(dt).
$$
The second term, through (\ref{roro}) is bounded by
$$
De^{qz}\exp(Be^{bz})\bb\int_{t_0}^1(1-t)^{\bb-1}dt=De^{qz}\exp(Be^{bz})(1-t_0)^{\bb}.
$$
Pick $t_0=t_0(z)>1-\varepsilon$ so that
\be \lbl{t0}\rho:=De^{qz}(1-t_0)^\bb<1.\ee 
In order to
establish (\ref{uppbeta}) we are to show that
$$e^{qz}\int_0^{t_0}\exp(Be^{bzt})\mu(dt) \le(1-\rho)\exp(Be^{bz}).$$
It follows from (\ref{t0}) that 
$$t_0=1-e^{-qz/\bb}(\rho/D)^{1/\bb},$$
and thus for sufficiently large $z$ we have that
$t_0>1-\varepsilon$. Hence, the left--hand side above, by
(\ref{roro}) again, is bounded by 
$$e^{qz}\exp(Be^{bzt_0})\mu(0,t_0)\le e^{qz}\exp(Be^{bzt_0})\left(1-\frac{d}{D}\rho e^{-qz}\right),
$$
and we want this to be less or equal than $(1-\rho)\exp(Be^{bz})$.
Divide both sides by $\exp(Be^{qz})$ so that the inequality to be
proved reads
$$e^{qz}\exp\left(Be^{bzt_0}-Be^{bz}\right)\left(1-\frac{d}{D}\rho e^{-qz}\right)\le 1-\rho.$$
We drop the factor $1-\frac{d}{D}\rho e^{-qz}$ on the left and
look at the exponent. It is 
$$qz+Be^{bz(1-e^{-qz/\bb}(\rho/D)^{1/\bb})}-Be^{bz}=qz+Be^{bz}\left(e^{-bze^{-qz/\bb}(\rho/D)^{1/\bb}}-1\right).
$$
Set $b:=q/\bb$. 
Since $\rho/D<1$ 
we
have $bze^{-qz/\bb}(\rho/D)^{1/\bb}=bze^{-bz}(\rho/D)^{1/\bb}<
bze^{-bz}\le e^{-1}<\ln 2$. Since $e^{-u}-1\le -u/2$ for
$0<u<\ln2$ we see that the expression above is bounded by
$$
qz-Bbz\rho^{1/\bb}e^{bz}e^{-bz}/2=
qz\left(1-\frac{B\rho^{1/\bb}}{2\bb}\right),$$
and it is clear that
$$
e^{qz}\exp\left(Be^{bzt_0}-Be^{bz}\right)
\le\exp\left(qz( 1-\frac{B\rho^{1/\bb}}{2\bb})\right),
$$
 can be made arbitrarily small by increasing $B$ if necessary.
In particular, we can ensure that it is less than
$1-\rho$
for all $z$ not too close to 0. Thus, (\ref{betatail}) is proved with $b=q/\bb$.

To get the matching lower bound  note that using again instead of $M$
the equivalent law ${\rm beta}(1,\bb)$\ with the cdf $F(t)=1-(1-t)^\bb$ we have
$$\nu_\delta=1-F(1-\delta)=\delta^\bb.$$ Thus, by (\ref{roro})
\begin{eqnarray*}\P(R\ge x)&\ge&\left(d\frac{cq}x\right)^{\bb\frac{\ln(1-c)}{\ln(1-cq/x)}}
=\exp\left(-\bb\frac{\ln(1-c)}{\ln(1-cq/x)}(\ln x-\ln(dcq))\right)\\
&=& \exp\left(\bb\frac{\ln(1-c)}{cq}(x\ln x)(1+o(1))\right).
\end{eqnarray*}
Hence, by letting $c\to 0_+$ we get that
$$\liminf_{x\to\infty}\frac{\ln\P(R\ge x)}{x\ln x}\ge-\frac\bb q.$$
\qed

\section{Generalized beta$(1,\beta)$ distributions}
\begin{figure}
\begin{center}
{\bf a}\epsfig{figure=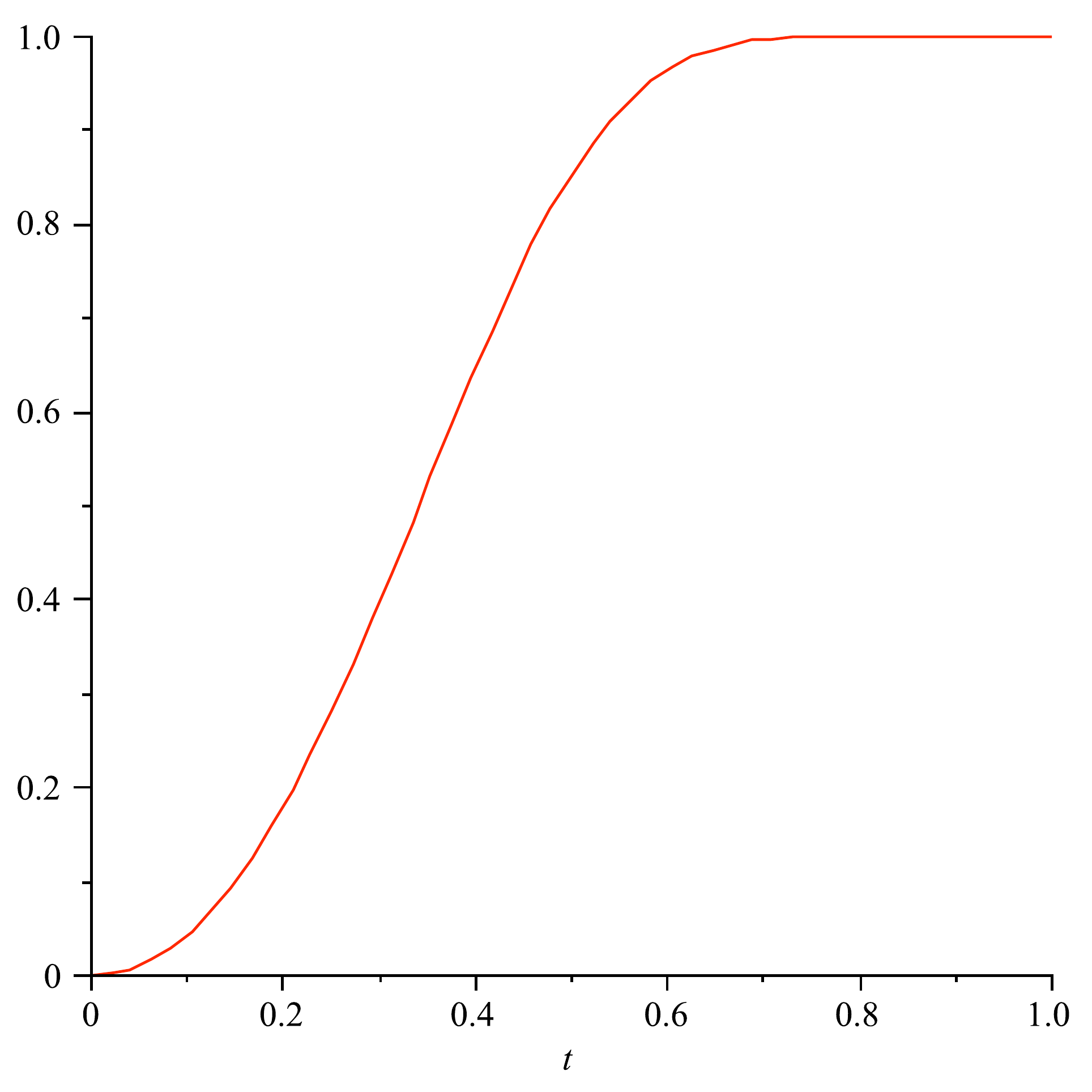, height=1.8in, width=2.0in, angle=0}\quad
{\bf b}\epsfig{figure=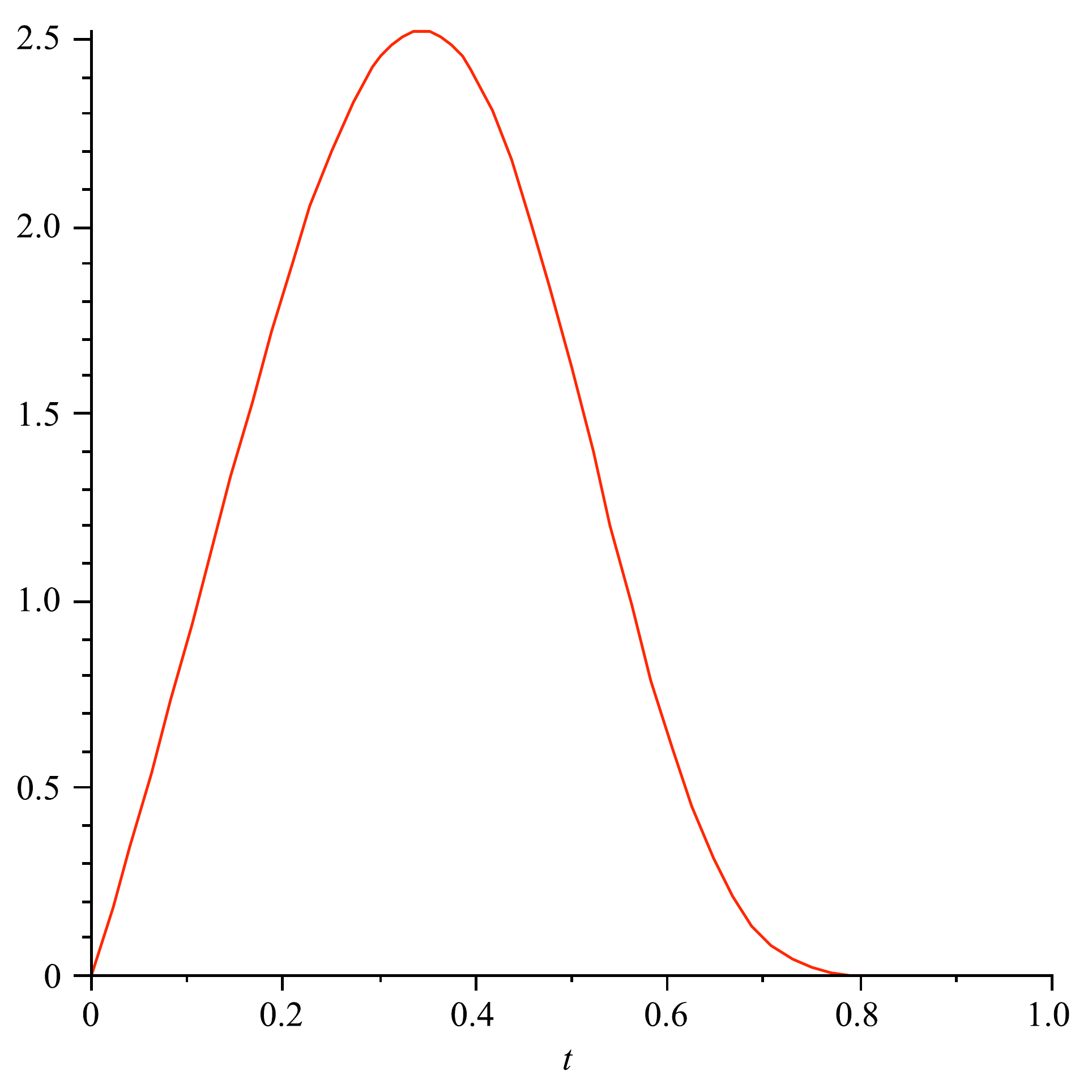, height=1.8in, width=2.0in, angle=0} \\
{\bf c}\epsfig{figure=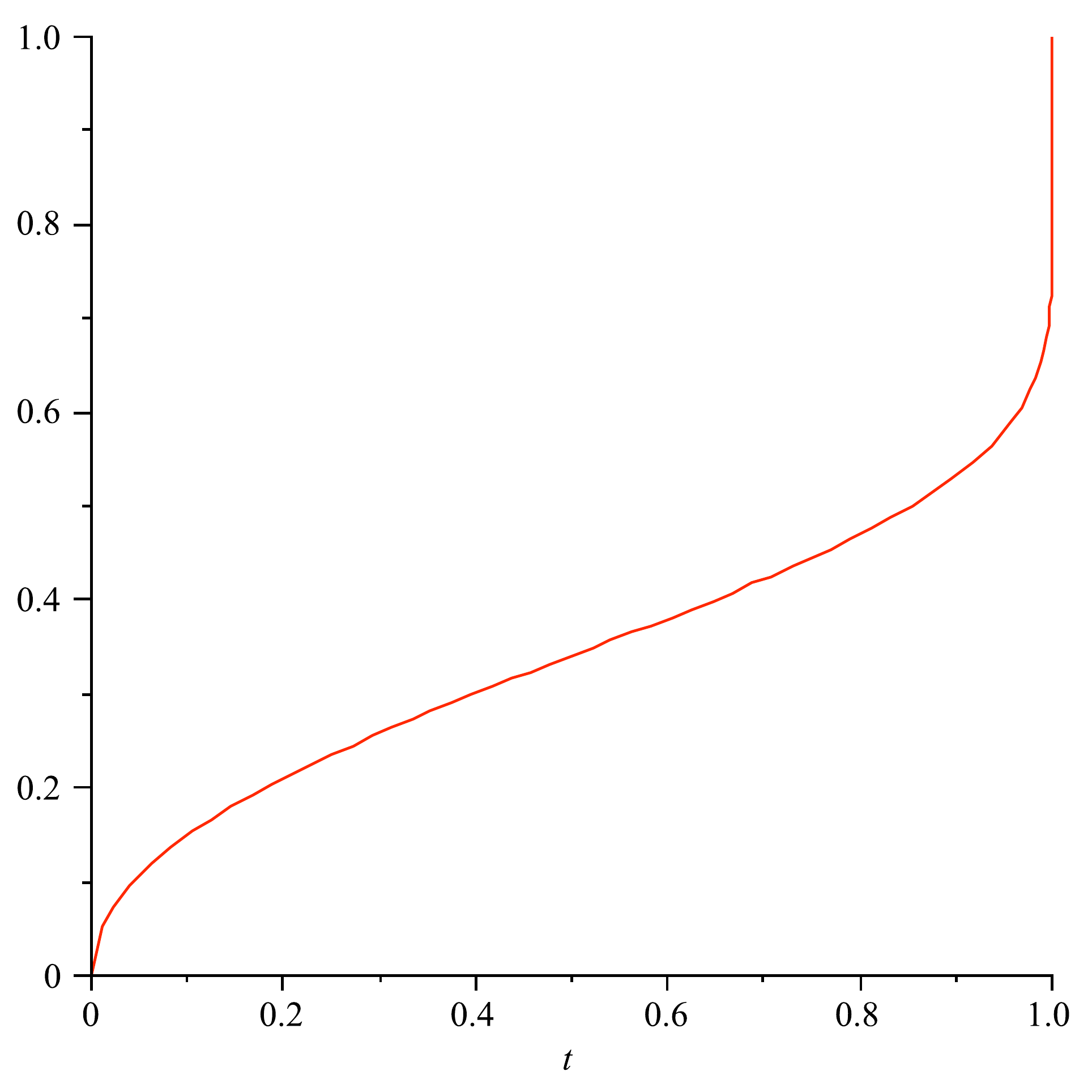, height=1.8in, width=2.0in, angle=0}\quad
{\bf d}\epsfig{figure=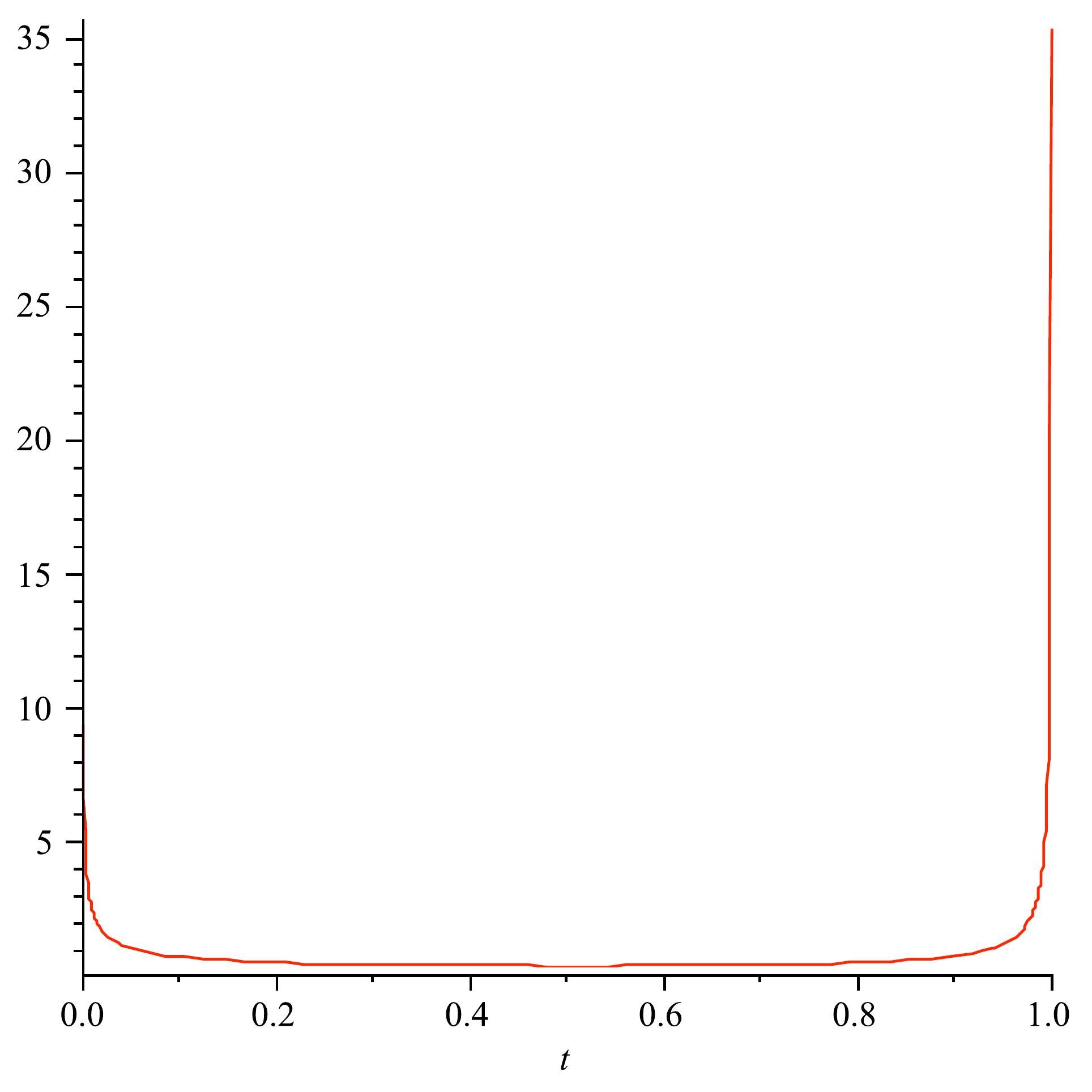, height=1.8in, width=2.0in, angle=0}
\end{center}
\caption{(a) The distribution $F_{4,2}$,\   (b) its density,\   (c) its  inverse $F_{0.5,0.5}$,\  and (d) its density.}
\label{fig:gb1}
\end{figure}

\begin{figure}
\begin{center}
{\bf a}\epsfig{figure=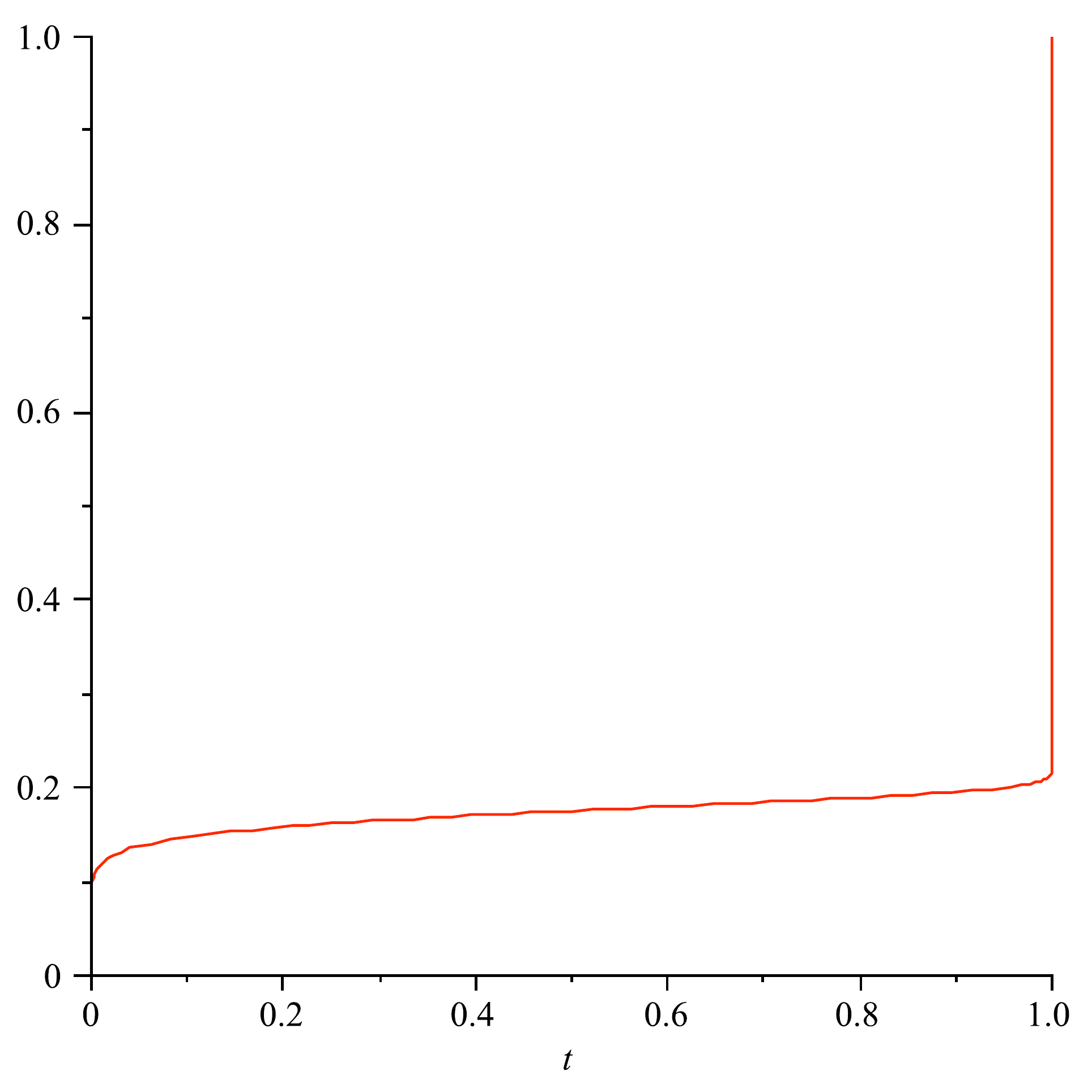, height=1.8in, width=2.0in, angle=0}\quad
{\bf b}\epsfig{figure=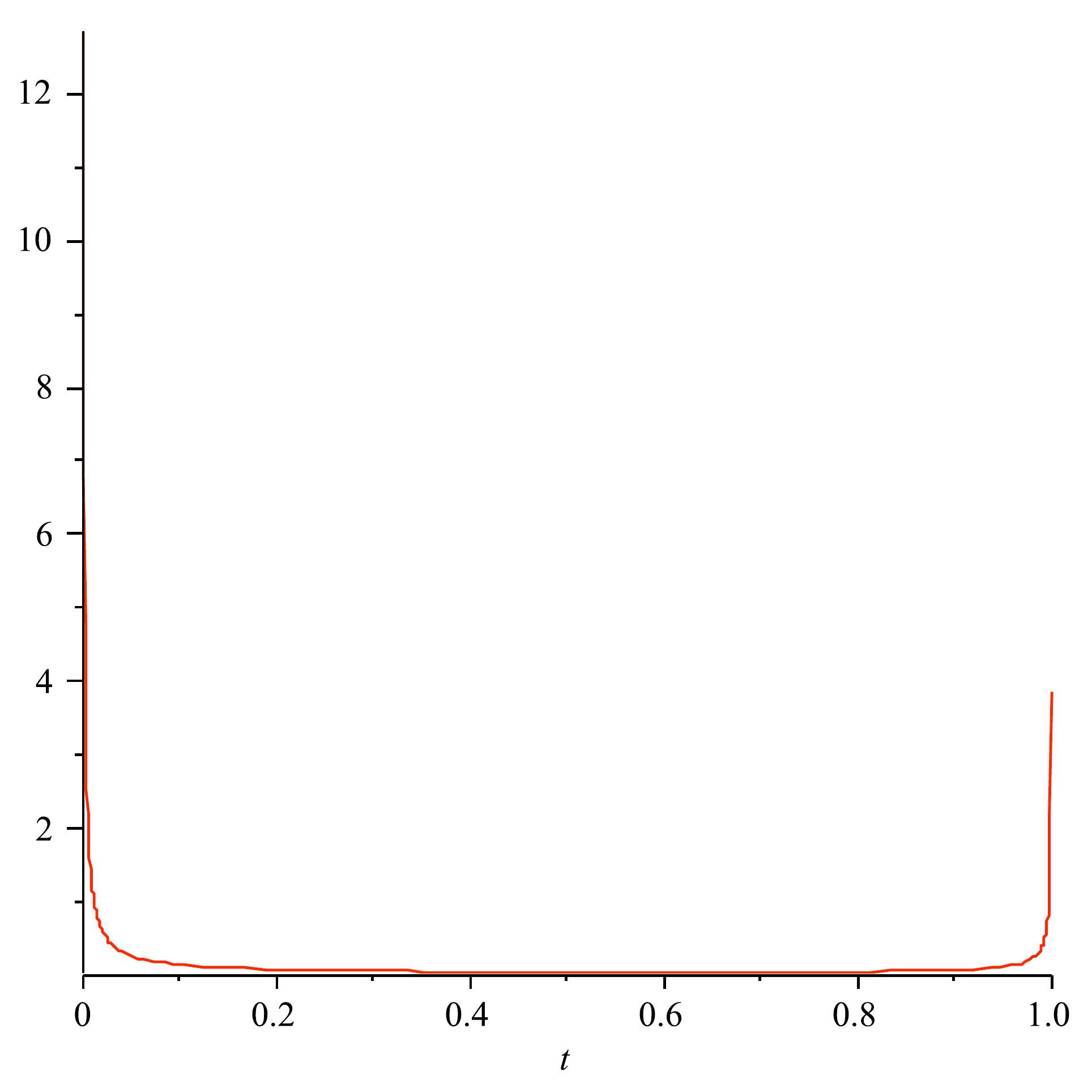, height=1.8in, width=2.0in, angle=0} \\
{\bf c}\epsfig{figure=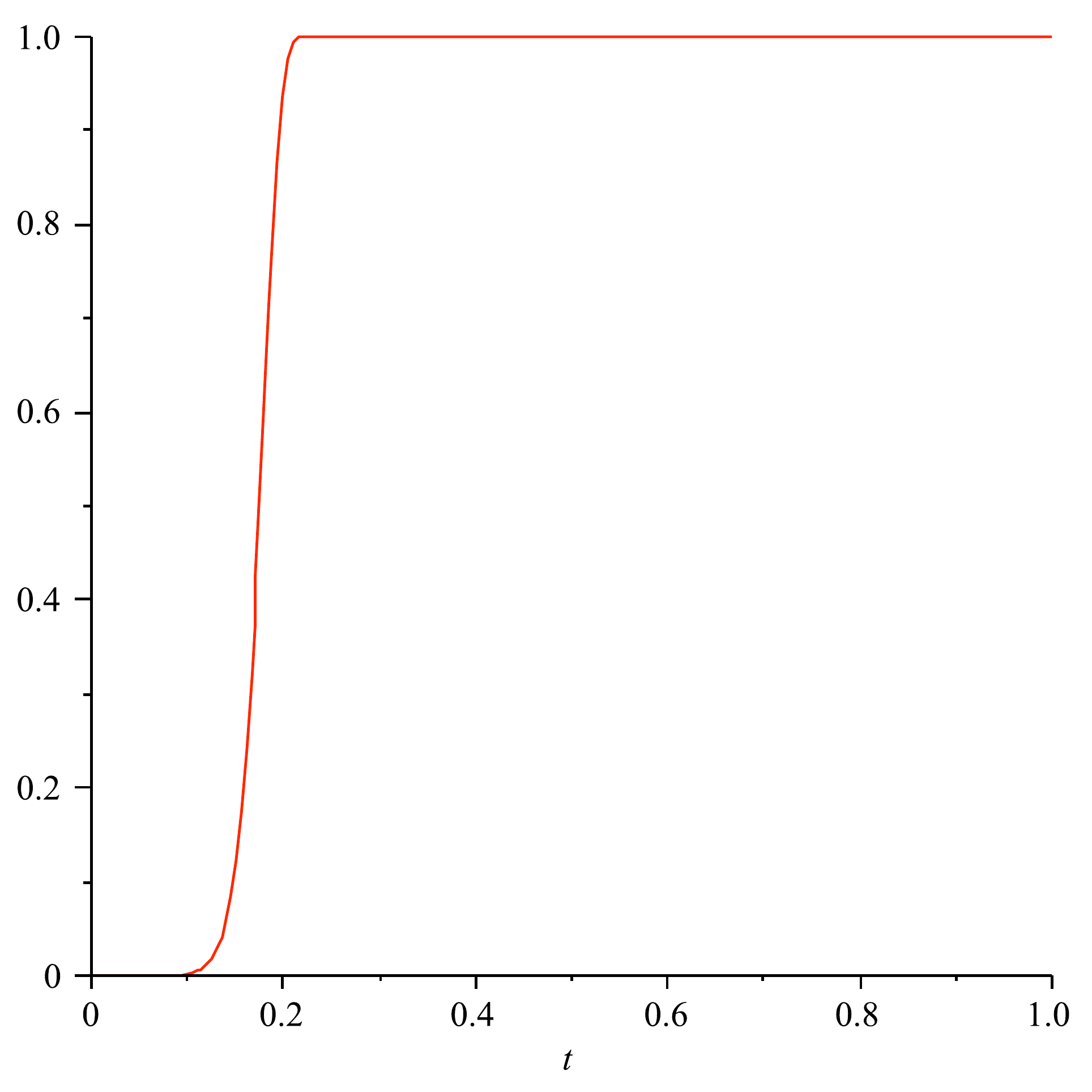, height=1.8in, width=2.0in, angle=0}\quad
{\bf d}\epsfig{figure=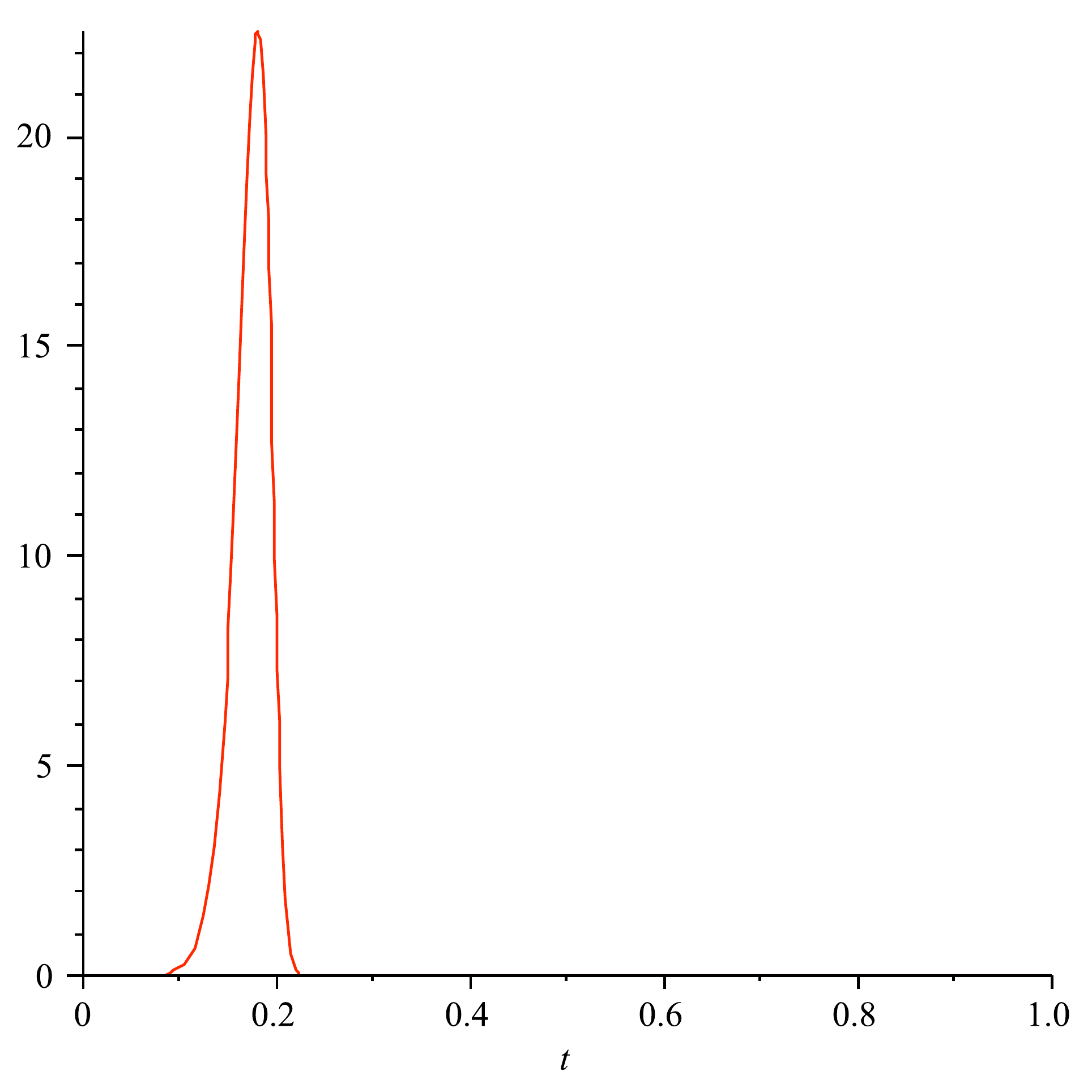, height=1.8in, width=2.0in, angle=0}
\end{center}
\caption{(a) The distribution $F_{0.2,0.1}$,\  (b) its density,\  (c) its inverse $F_{5^{10},10}$,  \    and (d) its  density.}
\label{fig:gb2}
\end{figure}

In this section we consider $M$'s whose distributions are
equivalent in the sense (\ref{roro}) to distribution function
given by
\be\lbl{gbeta}F(s)=F_{\bb,\eta}(s)=1-e^{-\bb(-\ln(1-s))^{\eta}},\quad0<s<1,\quad
\bb,\eta>0.\ee It is elementary to verify that $F_{\bb,\eta}$ is
indeed a distribution function which is strictly increasing on $(0,1)$.
Furthermore, $F_{\bb,1}$ is the distribution of a
$\mbox{beta}(1,\bb)$ random variable discussed in the previous
section. The family $F_{\bb,\eta}$ has the following property
$$F^{-1}_{\bb,\eta}=F_{\bb^{-1/\eta},\eta^{-1}},$$
as can be easily verified by a direct calculation.  Pictures of
a few such distributions with various parameters are given
in Figures~\ref{fig:gb1}--\ref{fig:gb2}.

For $R$ generated with $M$'s with distributions equivalent to  the
above distribution function  the following extension of
Theorem~\ref{thm:beta} holds
\begin{thm}
Let $(R_n)$ be given by (\ref{perp-n}) where $q>0$ and $M$ has the
distribution equivalent to the distribution function (\ref{gbeta})
for some $\bb,\eta>0$. Let $R$ be a limit in distribution of
$(R_n)$. Then
$$\lim_{x\to\infty}\frac{ \ln \P(R\ge x)}{x(\ln x)^{\eta}}= -\frac\bb q.$$
\end{thm}
\pf For the upper bound we will show that $R$ satisfies
Proposition~\ref{prop:upbdd} with $\Phi(z)=\exp(bz^{1/\eta})$ for
$b$'s in a certain range. For this $\Phi$  we have
\begin{eqnarray*}\Phi^*(x)\ge x\left(\left(\frac{\ln x}b\right)^\eta-B\right)\end{eqnarray*}
 which can be seen by using $\Phi_B^*(x)\ge xz_0-Be^{bz_0^{1/\eta}}$ with $z_0=b^{-\eta}(\ln x)^\eta$.
 It follows that
 \be\lbl{genbetalimsup}\limsup_{x\to\infty}\frac{\ln \P(R\ge x)}{x(\ln x)^\eta}\le -\frac1{b^\eta}.\ee
 To verify (\ref{iteration})
 we will use the same argument as before;  with $\Phi(z)=\exp(bz^{1/\eta})$ it  becomes
 $$e^{qz}\int_0^1\exp(Be^{b(zt)^{1/\eta}})\mu(dt)\le \exp(Be^{bz^{1/\eta}}),$$
 where $\mu$ is the distribution of the rv $M$ and $b$ and $B$ are positive constants.
Splitting the left--hand side, with $t_0>1-\varepsilon$ as before
we have
$$
\bb e^{qz}\int_0^{t_0}\exp(Be^{b(zt)^{1/\eta}})\mu(dt) +
e^{qz}\int_{t_0}^1\exp(Be^{b(zt)^{1/\eta}})\mu(dt).
$$
By (\ref{roro}) the second term is bounded by
$$
D e^{qz}\exp\left(Be^{bz^{1/\eta}}\right)(1-F_{\bb,\eta}(t_0))\;.
$$
Choose $t_0$ so that $\rho:=D
e^{qz}(1-F(t_0))<1$. Then
\begin{eqnarray*}t_0&=&F_{\bb,\eta}^{-1}(1-\rho e^{-qz}/D)\\ &=&F_{\bb^{-1/\eta},\eta^{-1}}(1-\rho e^{-qz}/D)=
1-\exp\left(-\bb^{-1/\eta}(-\ln(\rho e^{-qz}/D))^{1/\eta}\right)
\\&=&1-\exp\left(-\left(\frac{qz}\bb\right)^{1/\eta}\left(1-\frac{\ln(\rho/D)}{qz}\right)^{1/\eta}\right)
,
\end{eqnarray*}
and for $z$ sufficiently large it follows that
$t_0>1-\varepsilon$. Now we are to prove that
 $$e^{qz}\exp\left(Be^{bz^{1/\eta}t_0^{1/\eta}}\right)\mu(0,t_0)\le
(1-\rho)\exp\left(Be^{bz^{1/\eta}}\right)\;.
$$
By the first part of (\ref{roro}), it is enough to show that
\be\lbl{genbeta}
e^{qz}\exp\left(Be^{bz^{1/\eta}}\left(e^{-bz^{1/\eta}(1-t_0^{1/\eta})}
-1\right)
\right)\left(1-\frac{d\rho}{D}e^{-qz-Be^{bz^{1/\eta}}}\right) \le
1-\rho. \ee
We drop the last factor on the left--hand side as it is less that 1. For $t_0$ as
above $z^{1/\eta}(1-t_0^{1/\eta})$ is close to 0 for $z$
sufficiently large, so that using approximations $e^{-x}-1\sim-x$
and then $1-(1-x)^\eta\sim x/\eta$, both valid for $x$ close to 0
we see that the exponent on the left--hand side for $z$
sufficiently large,  is
\begin{eqnarray*}&&qz +Be^{bz^{1/\eta}}\left(e^{-bz^{1/\eta}(1-t_0^{1/\eta})}   -1\right)
\sim qz -Bbz^{1/\eta}e^{bz^{1/\eta}}(1-t_0^{1/\eta}) \\&&\quad
\sim qz -\frac{Bb}\eta z^{1/\eta}
\exp\left(z^{1/\eta}\left\{b-\left(\frac{q}\bb\right)^{1/\eta}\left(1-\frac{\ln\rho/D}{qz}\right)^{1/\eta}\right\}\right)\\&&\quad
\sim qz -\frac{Bb}\eta z^{1/\eta}
\exp\left(z^{1/\eta}\left\{b-(q/\bb)^{1/\eta}\right\}\right).
\end{eqnarray*}
For $b>(q/\bb)^{1/\eta}$ the second term grows faster than
linearly in $z$, so that as long as $z$ is not too close to 0 it
can be made arbitrarily larger than $qz$. Thus, (\ref{genbeta})
follows. Furthermore, letting $b\to(q/\bb)^{1/\eta}_+$ in
(\ref{genbetalimsup}) we obtain that \be
\lbl{lsup}\limsup_{x\to\infty}\frac{ \ln \P(R\ge x)}{x(\ln
x)^{\eta}}\le -\frac\bb q. \ee To get a lower bound note that,
using instead of the distribution of $M$ the equivalent cdf
$F_{\bb,\eta}$, on noting that
$$1-F_{\bb,\eta}(1-cq/x)=\exp\left(-\bb\left(-\ln(cq/x)\right)^\eta\right)=
\exp\left(-\bb(\ln x-\ln(cq))^\eta\right)
$$
we get for large $x$
\begin{eqnarray*}
\P(R\ge
x)&\ge&\left(d(1-F_{\bb,\eta}(1-cq/x))\right)^{\frac{\ln(1-c)}{\ln(1-cq/x)}}\\ &=&\exp\left(-\frac{\ln(1-c)}{\ln(1-cq/x)}\bb[(\ln
x-\ln (cq))^\eta+\ln(d)]\right) \\ &=&
\exp\left(\frac{\bb\ln(1-c)}{cq}x(\ln x)^\eta(1-o(1))\right).
\end{eqnarray*}
Upon letting $c\to0_+$ it implies that
$$
\liminf_{x\to\infty}\frac{ \ln \P(R\ge x)}{x(\ln x)^{\eta}}\ge
-\frac\bb q\;.
$$
Combining this with  (\ref{lsup}) completes the proof.\qed

\section{Weilbull--like tails}
In this section we explicitly construct  $M$'s that will lead to a rather different tail behavior of $R$ than discussed in the previous sections. As we will see a much more rapid variability of $M$ near 1 is needed to obtain a lighter tail behavior of $R$. More specifically, we prove the following theorem.
\begin{thm}
Let $1<r<\infty$. Let the distribution of $M$ be (\ref{roro})
equivalent to the distribution $\nu$  with the density 
\be\lbl{density} f_\nu(t)\ \propto\  
t^{1/(r-1)}e^{-\frac
{1}{(1-t^{r/(r-1)})^{r-1}}}I_{(0,1)}(t).
\ee Then for the perpetuity
$R$ given by (\ref{perp-it}) with $Q\equiv q$ there are constants
$c_1,c_2$ such that
$$
-\infty<c_1\le\liminf_{x\to\infty}\frac{\ln \P(R\ge x)}{(x/q)^r} \le
 \limsup_{x\to\infty}\frac{\ln \P(R\ge x)}{(x/q)^r}\le c_2<0.
$$
\end{thm}
\pf
For $1<r<\infty$ let $r^*$ be given by
$$
\frac1r+\frac1{r^*}=1.
$$
The role of $r$ and $r^*$ are symmetric and for notational
convenience we will prove the above inequalities for $r^*=r/(r-1)$ rather
than $r$. That is to say, we will show that if $M$ is equivalent to a random variable whose density is
proportional to 
\be\lbl{densityr*}
t^{1/(r^*-1)}e^{-\frac
{1}{(1-t^{r^*/(r^*-1)})^{r^*-1}}}I_{(0,1)}(t)
=
t^{r-1}e^{-\frac
{1}{(1-t^{r})^{1/(r-1)}}}I_{(0,1)}(t),
\ee 
then the resulting perpetuity $R$ satisfies
 $$
-\infty<c_1\le\liminf_{x\to\infty}\frac{\ln \P(R\ge x)}{(x/q)^{r^*}} \le
 \limsup_{x\to\infty}\frac{\ln \P(R\ge x)}{(x/q)^{r^*}}\le c_2<0,
$$
for some constants 
$c_1,c_2$.

Suppose we prove that for $M$ the condition
(\ref{mgf-n}) holds for all $n\ge 1$ with $\Phi(z)=z^r$ and some
$B>0$. Then by elementary calculation
$\Phi^*(x)=\frac{x^{r^*}}{r^*(Br)^{1/(r-1)}}$, so that,
\be\lbl{uppw}\P(R\ge
x)\le\exp\left(-\frac{x^{r^*}}{r^*(Br)^{1/(r-1)}}\right),\ee
and this would give the claimed behavior of the logarithm of the tail probability of $R$. 

To establish (\ref{mgf-n}) via inductive argument, we need to verify that (\ref{iteration}) holds in the present situation, that is, 
we want to show that for $z$ sufficiently large
$$
e^{qz}\int_0^1\:e^{B(zt)^r}\:\mu(dt)\le e^{Bz^r}\;.
$$
Take $\varepsilon>0$ given by (\ref{roro}) where $\nu$ has density of the form (\ref{densityr*}) and consider $\delta\in(0,\varepsilon)$. Then
the left hand side of the above inequality is less than
$$
e^{qz}e^{Bz^r(1-\delta)^r}+e^{qz}\int_{1-\delta}^1\:e^{B(zt)^r}\:\mu(dt)\le e^{qz}e^{Bz^r(1-\delta)^r}+ De^{qz}\int_{1-\delta}^1\:e^{B(zt)^r}\:\nu(dt)\;.
$$
Consequently, we have to show that
\be\lbl{nu:bdd}e^{qz-Bz^r(1-(1-\delta)^r)}+De^{qz-Bz^r}\int_{1-\delta}^1\:e^{B(zt)^r}\:f_\nu(t)dt\le 1\;.\ee
Note that because $r>1$ and $0<\delta<1$, the first term can be made arbitrarily small for $z\ge z_0$ sufficiently large. We thus  concentrate on the second term. 
The following argument will not only complete justification of (\ref{nu:bdd}) but will also indicate how one would be led to a reasonable choice of $f_\nu$ if it were unknown.  We would want to construct a density $f_\nu$ on
$(0,1)$ for which (\ref{nu:bdd}) holds. To
this end suppose for now that the density $f_\nu$  were of the form
$$f_\nu(t)=rt^{r-1}g(t^r).$$
Upon changing variables to $s=t^r$  the second term in (\ref{nu:bdd}) becomes
$$De^{qz-Bz^r}\int_{(1-\delta)^r}^1e^{Bz^rs}g(s)ds=De^{qz}\int_{(1-\delta)^r}^1e^{-Bz^r(1-s)}g(s)ds.$$
Setting $w=1-s$ gives
\be\lbl{lastsub}De^{qz}\int_0^{1-(1-\delta)^r}e^{-Bz^rw}g(1-w)dw
.\ee
We now let
$$
g(1-w):=Ke^{-1/w^\gg},$$
where  $\gg$ is to be chosen momentarily and $K=K(\gg)$ is set so that
$$K^{-1}=\int_0^1e^{-1/w^\gg}dw.$$
Then  (\ref{lastsub}) becomes
\be\lbl{int-ineq}KDe^{qz}\int_0^{1-(1-\delta)^r}e^{-Bz^rw}e^{-1/w^\gg}dw.
\ee
The integrand  is
$$\exp\left(-(Bz^rw+\frac1{w^\gg})\right).$$
Since the function
$$w\to Bz^rw+\frac1{w^\gg},$$
has a minimum at $(\gg/(Bz^r))^{1/(\gg+1)}$ whose value is
$$
(Bz^r)^{\frac\gg{\gg+1}}(\gg^{\frac1{\gg+1}}+\gg^{-\frac\gg{\gg+1}})=
B^{\frac\gg{\gg+1}}z^{r\frac\gg{\gg+1}}\frac{\gg+1}{\gg^{\gg/(\gg+1)}},
$$
the quantity (\ref{int-ineq}) is no more than
$$KD\exp\left(zq-z^{r\frac\gg{\gg+1}}B^{\frac\gg{\gg+1}}\frac{\gg+1}{\gg^{\gg/(\gg+1)}}\right),$$
which upon setting
$$r\frac\gg{\gg+1}=1\quad\mbox{i.e.}\quad
\gg=\frac1{r-1},$$
becomes
$$KD\exp\left\{z\left(q-B^{1/r}\frac r{(r-1)^{(r-1)/r}}\right)\right\}.
$$
It is now clear that if
\be\lbl{B}B=A^r\left(\frac qr\right)^r(r-1)^{r-1},\ee
where $A>1$ might depend on $r$, then 
$q-B^{1/r}r/(r-1)^{\frac{r-1}r}=q(1-A)<0$. 
Therefore, for $z\ge z_0$ we obtain further
$$KD\exp\left\{z\left(q-B^{1/r}\frac r{(r-1)^{(r-1)/r}}\right)\right\}\le KDe^{-z_0q(A-1)}.
$$
Thus we conclude that for $z\ge z_0$ the left--hand side of  (\ref{nu:bdd}) is bounded by 
$$e^{-z_0(B(1-(1-\delta)^r)z_0^{r-1}-q)}+KDe^{-z_0q(A-1)}.$$
Since the value of this expression can be made  smaller than 1
by choosing $z_0$ sufficiently large, (\ref{nu:bdd})
follows.

Reversing the steps, we obtain the expression for the density  
given in (\ref{densityr*})
with the normalizing constant $K_r$  given by
$$K_r^{-1}=\frac1r\int_0^1\exp\left(-\frac1{v^{1/(r-1)}}\right)dv.$$
\begin{figure}
\begin{center}
{\bf a}\epsfig{figure=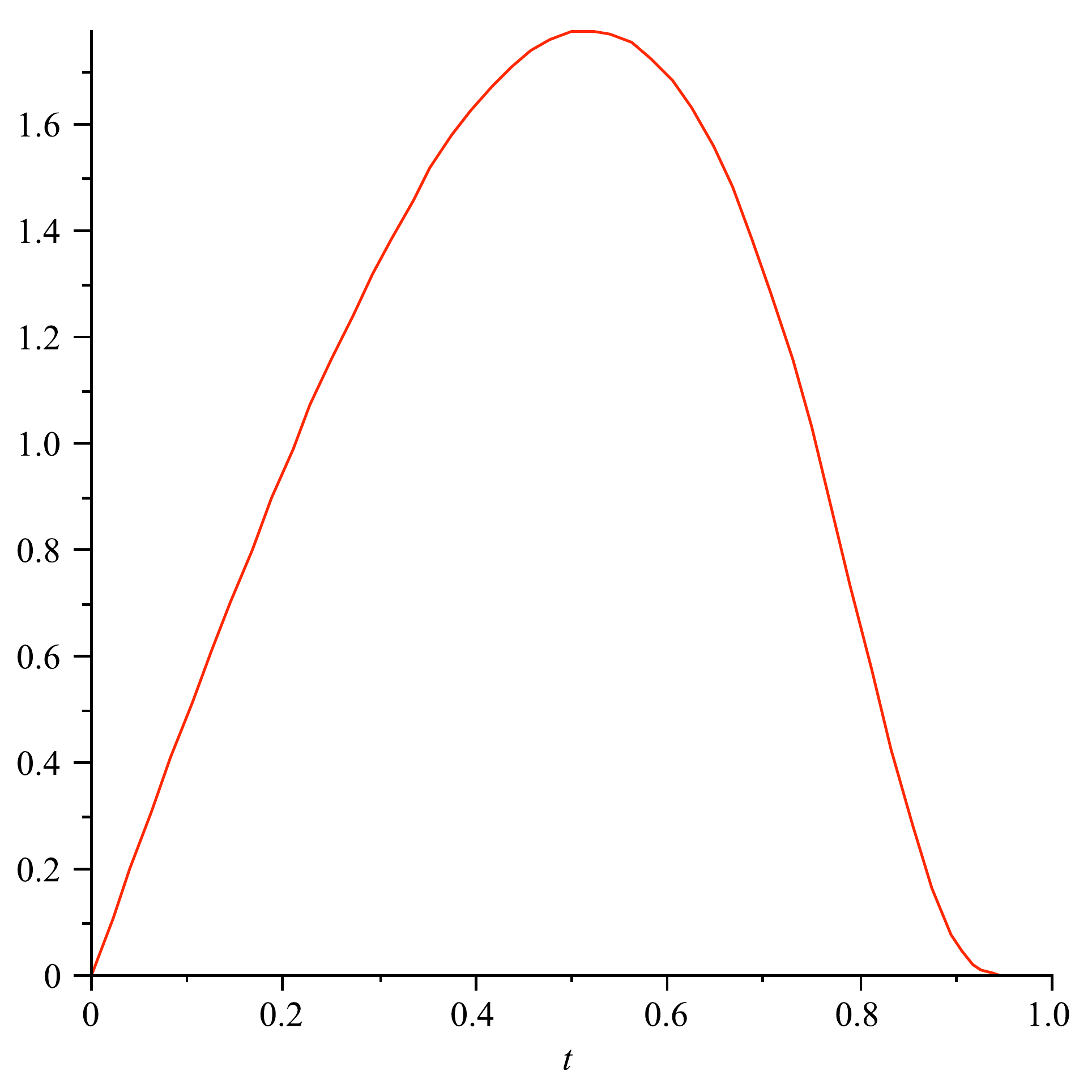, height=1.6in, width=2.0in, angle=0}\quad
{\bf b}\epsfig{figure=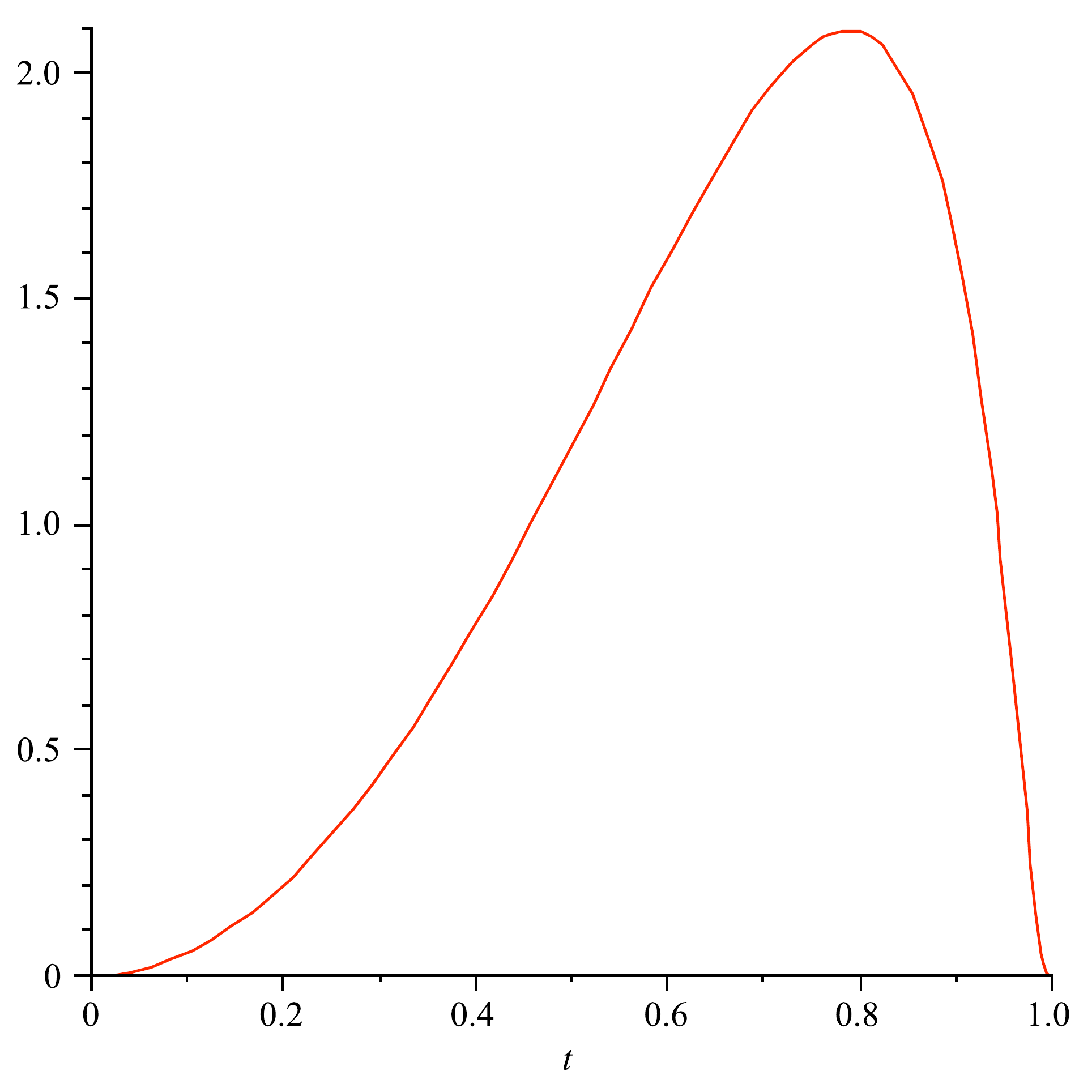, height=1.6in, width=2.0in, angle=0}
\end{center}
\caption{The density (\ref{densityr*}) for (a) $r=2$ \   and (b) $r=3$.
}
\label{fig:weil23}
\end{figure}

\begin{figure}
\begin{center}
{\bf a}\epsfig{figure=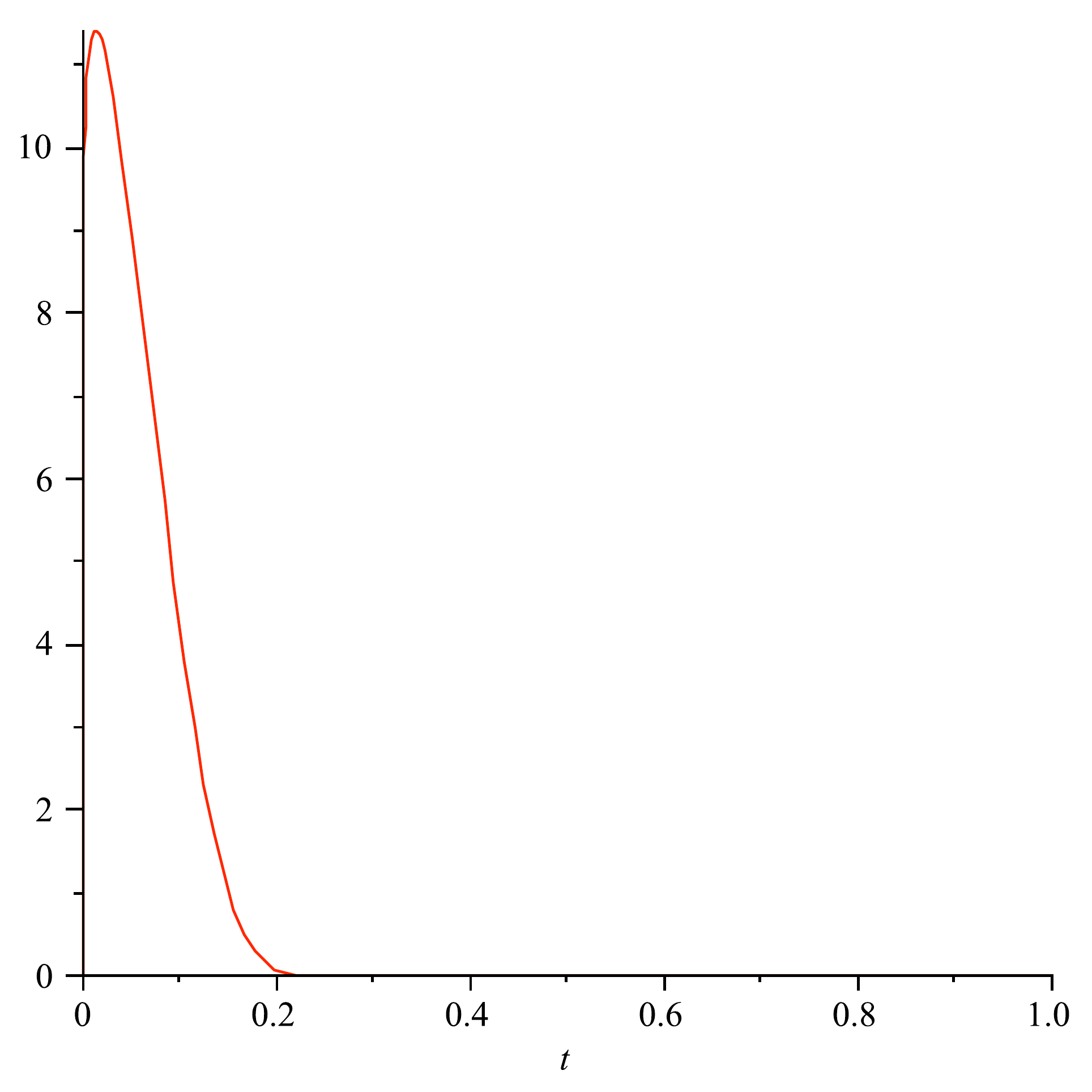, height=1.6in, width=2.0in, angle=0}\quad
{\bf b}\epsfig{figure=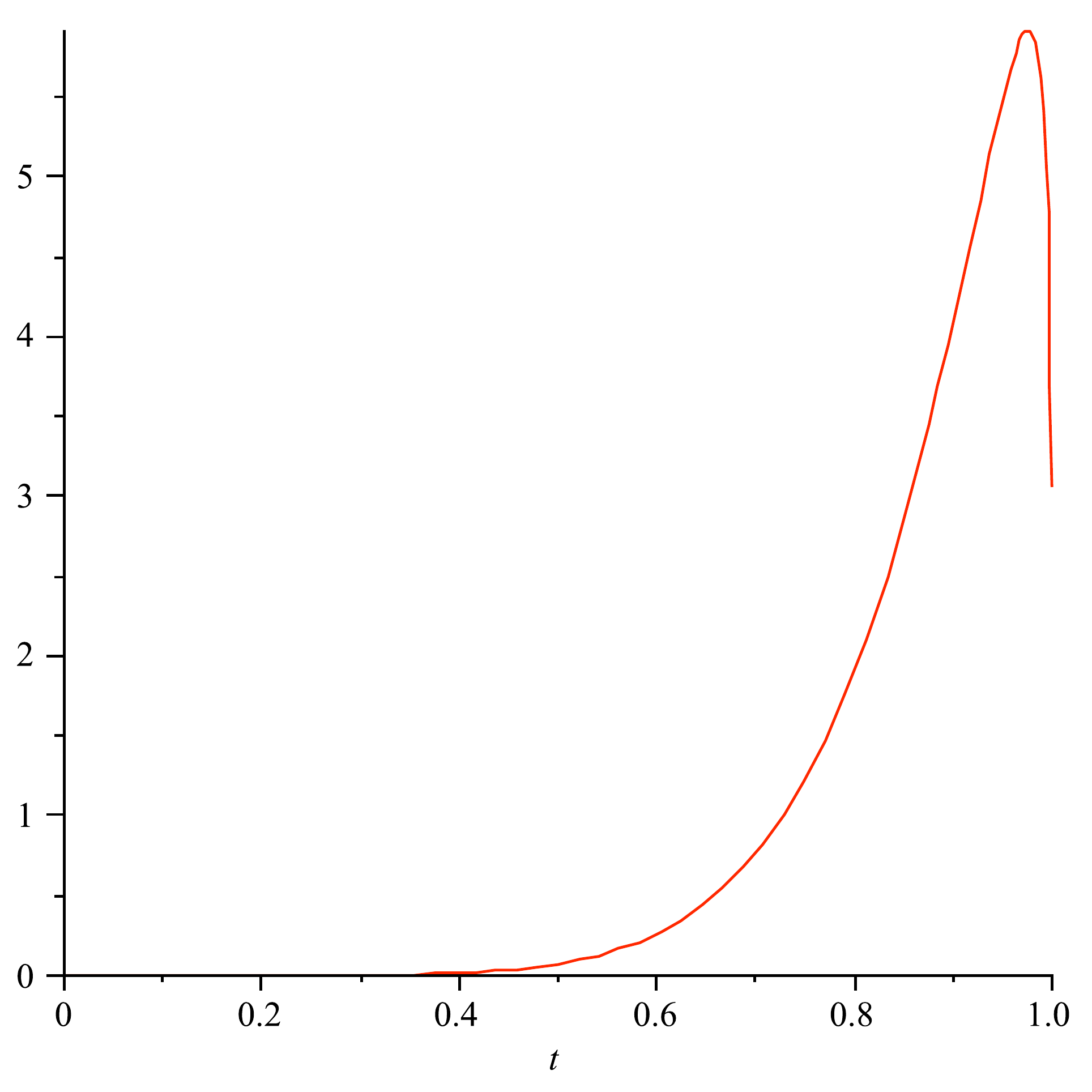, height=1.6in, width=2.0in, angle=0}
\end{center}
\caption{The density (\ref{densityr*}) for (a) $r=1.1$ \   and (b) $r=8$.
}
\label{fig:weil118}
\end{figure}

Finally, putting the value of $B$ given in (\ref{B}) into (\ref{uppw}) we obtain
$$\P(R\ge x)\le\exp\left(-\left(\frac xq\right)^{r^*}\frac1{A^{r/(r-1)}}\right),$$
which implies that
$$\limsup_{x\to\infty}\frac{\ln\P(R\ge x)}{(x/q)^{r^*}}\le -A^{-r/(r-1)}.$$

To get  a lower bound for $\P(R\ge x)$ we choose $\delta\in(0,\varepsilon)$ as in (\ref{roro}). Then, upon passing to the equivalent measure with density proportional to (\ref{densityr*}) we have
$$p_\delta \ge d K_r\int_{1-\delta}^1t^{r-1}\exp\left(-\frac1{(1-t^r)^{1/(r-1)}}\right)dt.$$
Changing variables to $v=(1-t^r)^{-1/(r-1)}$ yields 
$$p_\delta\ge K\int_{(1-(1-\delta)^r)^{-1/(r-1)}}^\infty\frac{e^{-v}}{v^r}dv,$$
for some constant $K$ whose value is irrelevant. 
Since for large $v_0$, $\int_{v_0}^\infty\frac{e^{-v}}{v^r}dv$ is comparable to
$e^{-v_0}/v_0^r$ we  get, up to an unimportant constant
$$(1-(1-\delta)^r)^{r/(r-1)}
\exp\left(-\frac1{(1-(1-\delta)^r)^{1/(r-1)}}\right)
$$
as the lower bound for $p_\delta$. Hence, up to unimportant additive terms
$$\ln p_\delta\ge\frac r{r-1}\ln\left(1-(1-\delta)^r\right)-\frac1{(1-(1-\delta)^r)^{1/(r-1)}}\sim -\frac1{(1-(1-\delta)^r)^{1/(r-1)}},$$
as the second term above is of dominant order for $\delta\to0$.
For small $\delta$ we have
$$1-\left(1-\delta\right)^r=1-\exp\left(r\ln(1-\delta)\right)\sim-r\ln(1-\delta),
$$
so that upon replacing $\delta$ by $cq/x$ we get that, asymptotically
$$\ln p_{cq/x}\ge -\frac1{(-r\ln(1-cq/x))^{1/(r-1)}}\sim-\left(\frac x{cqr}\right)^{\frac1{r-1}}.$$
Combining this with (\ref{lbdd-x}) we get that, asymptotically,
\begin{eqnarray*}\ln \P(R\ge x)&\ge&\frac{\ln(1-c)}{\ln(1-cq/x)}  \left(-\left(\frac x{cqr}\right)^{1/(r-1)}\right)
\sim\frac{x\ln(1-c)}{cq} \left(\frac x{cqr}\right)^{1/(r-1)}  \
\\&=&
\left(\frac xq\right)^{r^*}\cdot\frac{\ln(1-c)}{(cr^{1/r})^{r^*}}.
\end{eqnarray*}
It follows that
$$\liminf_{x\to\infty}\frac{\ln\P(R\ge x)}{(x/q)^{r^*}}\ge  \frac C{r^{1/(r-1)}},\quad\mbox{where}\quad C=\frac{\ln(1-c)}{c^{r^*}}<0. $$\qed
\noindent{\bf Remarks:}

\noindent (i) The  maximal  value of $C/r^{1/(r-1)}$ is obtained by setting $c=c_0$ where $c_0 $ is  the unique solution of  the equation
$$ \frac1{1-c}+ r^*\cdot\frac{\ln(1-c)}c=0.$$
The uniqueness of the solution is elementary as the function
$$h(c):=\frac{\ln(1-c)}{c^{r^*}}$$
approaches $-\infty$ as $c\to0_+$ or $c\to 1_-$ and
$$h'(c)=-c^{-r^*}\left(\frac1{1-c}+r^*\cdot\frac{\ln(1-c)}c\right).$$
The expression in the parentheses, upon letting $y=1/(1-c)$, $y>1$, becomes
$$y-r^*\cdot\frac{\ln y}{(y-1)/y}=y\left(1-r^*\cdot \frac{\ln y}{y-1}\right).$$
Since $\frac{\ln y}{y-1}$ is decreasing for $y>1$, approaches $1$ as $y\to 1_+$ and $0$ as $y\to\infty$ we see that $h'(c)$ has exactly one sign change (from positive to negative) on $(0,1)$ and that this change occurs at $c_0$ such that
$$\frac1{1-c_0}+r^*\cdot\frac{\ln(1-c_0)}{c_0}=0.$$
While the above equation does not have in general the closed form solution for $c_0$ as a function of $r$ (or $r^*$), the asymptotic behavior of the constant $C/r^{1/(r-1)}$ as $r$ goes to  0 or $\infty$ can  be traced down. Since
$$r^*=-\frac{c_0}{(1-c_0)\ln(1-c_0)},$$
as $r\to\infty$ (and thus $r^*\to 1_+$) we must have $c_0\to0_+$ at the rate
$1-c_0\sim 1/r^*$. But then $c_0\sim1-1/r^*=1/r$ and thus
$$\frac{\ln(1-c_0)}{r^{1/(r-1)}c_0^{r^*}}\sim\frac{\ln(1-1/r)}{r^{1/(r-1)}(1/r)^{r/(r-1)}}=
r\ln(1-1/r)\to-1,\quad \mbox{as\quad $r\to\infty$}.
$$
Similarly, if $r\to1_+$  then $c_0\to1_-$ in such a way that $1-c_0\sim1/(r^*\ln r^*)$.
Then
$$\frac{\ln(1-c_0)}{r^{1/(r-1)}c_0^{r^*}}\sim\frac{-\ln(r^*\ln r^*)}{r^{1/(r-1)}(1-1/(r^*\ln r^*))^{r^*}}
\sim\frac{-\ln(r^*\ln r^*)}e
$$
since, as $r\to1_+$
$$
r^{1/(r-1)}=\left(1+\frac1{1/(r-1)}\right)^{1/(r-1)}\to e,\quad \mbox{and}
\quad \left(1-\frac1{r^*\ln r^*}\right)^{r^*}\to 1.
$$

\noindent(ii) It might appear from the argument that the form of density (\ref{density}) was just guessed. While it is true that originally this was the case, there is a heuristic argument which would suggest the same choice. We will explain this heuristics in the next section on a different example, but we would like to mention that following it in the present situation would essentially lead to  density given by (\ref{density}).
\section{Further example} In this section we present one more example  of perpetuity that will have extremely thin tails. Specifically, we will show
\begin{pr}
There exist densities $f_M$ for which the  perpetuity defined by (\ref{perp-n}) satisfies:
\be\lbl{upbdd-j}\forall\ B>q\quad\limsup_{x\to\infty}\frac{\ln \P(R\ge x)}{B\exp(x/B)}\le-\frac1e,
\ee
and
\be\lbl{lobdd-j}
\forall\ B<q\quad\liminf_{x\to\infty}\frac{\ln \P(R\ge x)}{B\exp(x/B)}\ge\frac{\ln(1-B/q)}B.
\ee
\end{pr}
\pf We consider the case
$\Phi(z)=z\ln z$ and we will  show that Proposition~\ref{prop:upbdd} holds for all $B>q$.   It will then follow that for all such $B$
\be\lbl{upextr}\P(R\ge x)\le\exp(-B\exp(\frac xB-1)),\ee
which will  imply (\ref{upbdd-j}).
We will then construct a density of $M$ which, on one hand will  guarantee (\ref{upextr}) and, on the other hand, ensure that $p_\delta$ is sufficiently large so that the argument based on Proposition~\ref{prop:lbdd} will give (\ref{lobdd-j}).  

\begin{figure}
\begin{center}
{\bf a}\epsfig{figure=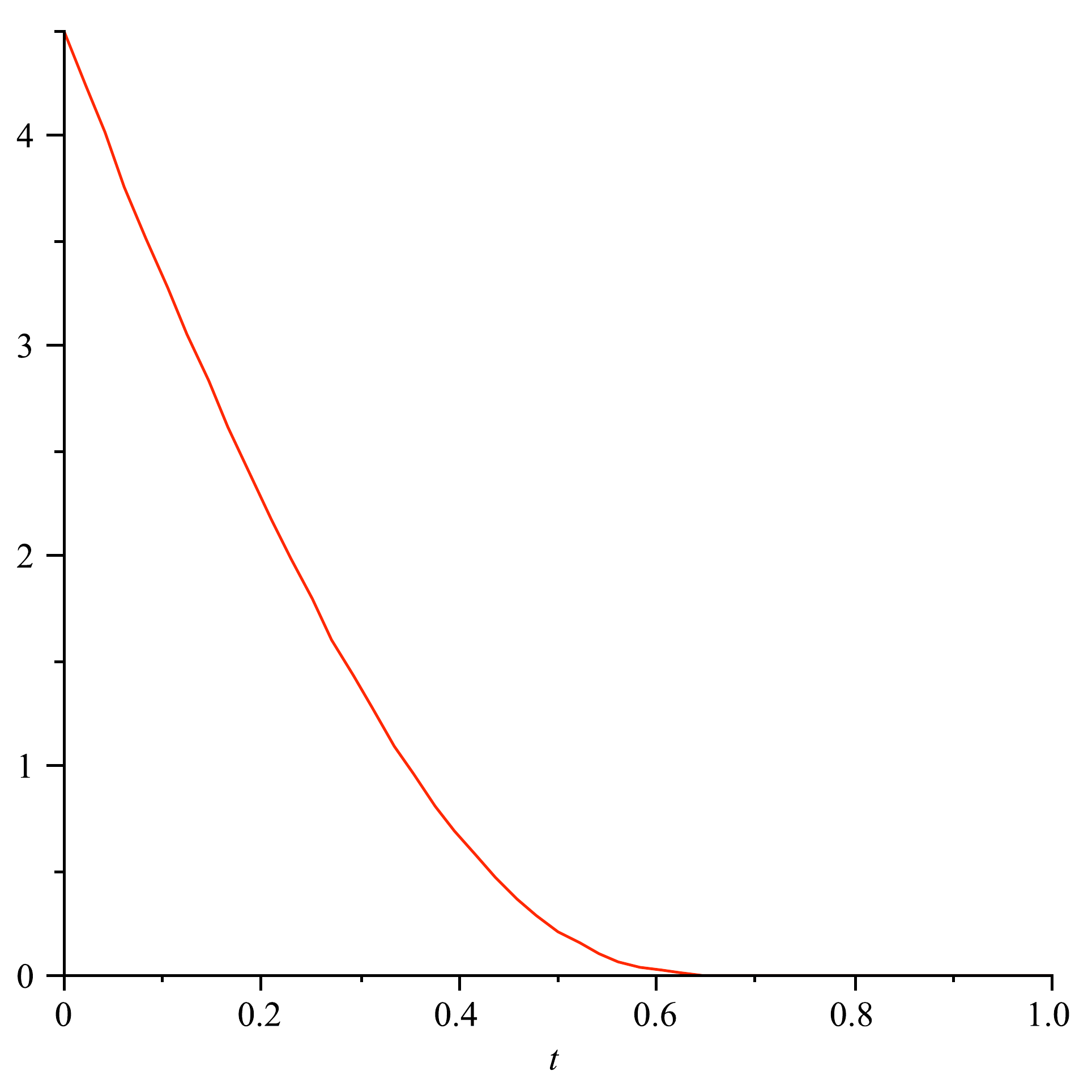, height=1.6in, width=2.0in, angle=0}\quad
{\bf b}\epsfig{figure=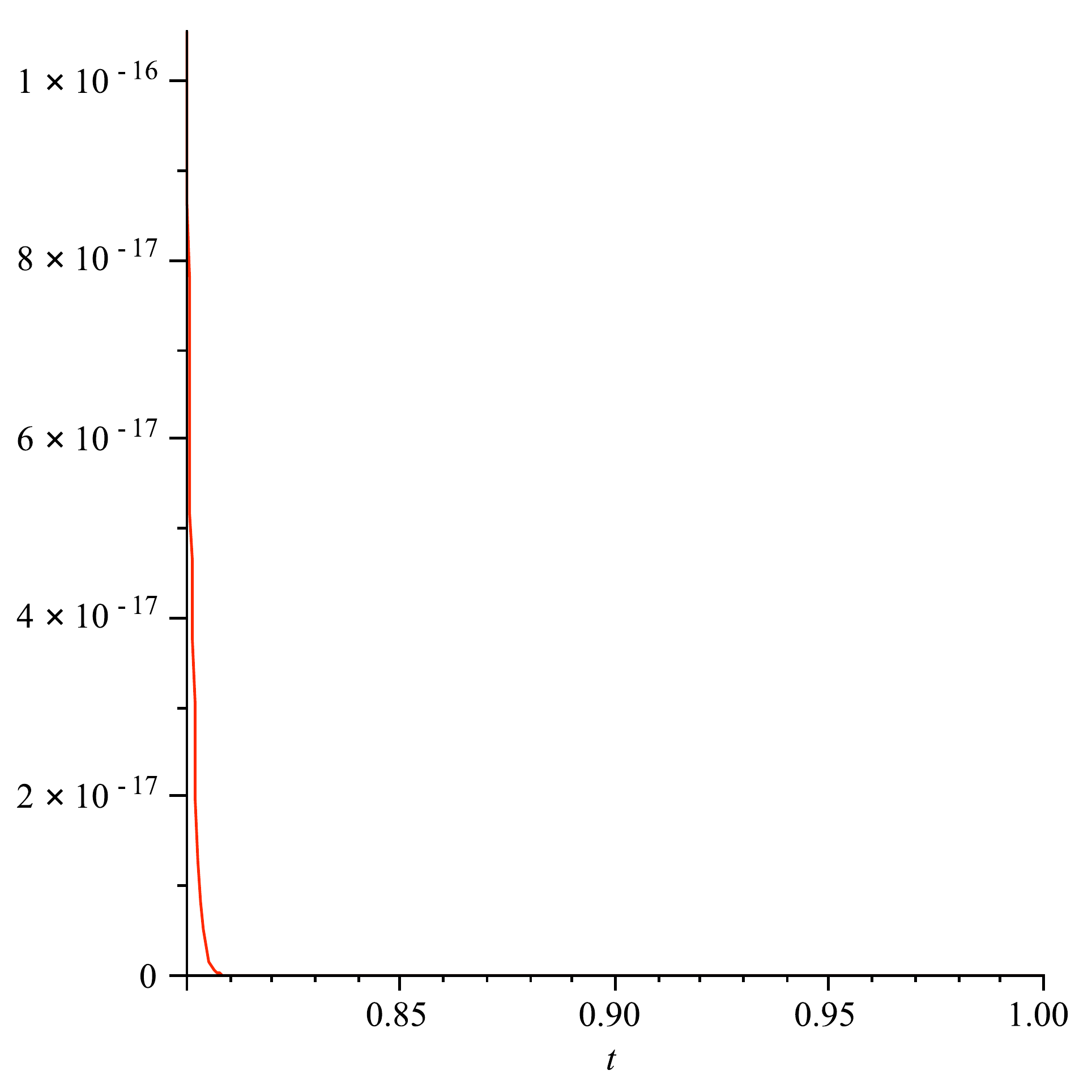, height=1.6in, width=2.0in, angle=0}
\end{center}
\caption{(a) The density $f_M$\   and (b) its  detail  closer to 1.
}
\label{fig:jdens}
\end{figure}

To carry out the details of that plan we are to construct a density $f_M$ for which
$$e^{qz}\int_0^1e^{Bzt\ln(zt)}f_M(t)dt\le e^{Bz\ln z}.$$
This is equivalent to
$$e^{qz}\int_0^1e^{-B(1-t)z\ln z}t^{Btz}f_M(t)dt\le 1,$$
and it is enough to construct an $f_M$ for which
$$e^{qz}\int_0^1e^{-B(1-t)z\ln z}f_M(t)dt=
e^{qz}\int_0^1e^{-Btz\ln z}f_M(1-t)dt
\le 1.$$
We now set $f_M(1-t)=K\exp(-h(t)),$ where $h$ is a non--negative function and $K^{-1}=\int_0^1\exp(-h(t))dt$. The inequality to be established becomes
\be\lbl{jacekineq}e^{qz}\int_0^1e^{-Btz\ln z-h(t)}dt\le \int_0^1e^{-h(t)}dt.
\ee
One is guided to a reasonable choice of $h$ by the following heuristics.  Suppose $h$ is  differentiable and chosen so that
\be\lbl{et}Btz\ln z+h(t)\ee
  is minimized at its critical point $t=t_z\in(0,1)$  which thus satisfies
\be\lbl{tofz}Bz\ln z+h'(t_z)=0.\ee
Then the left--hand side of (\ref{jacekineq}) is no more than
  $$\exp\left(qz-Bt_zz\ln z-h(t_z)\right)
  \le
\exp\left(z(q-Bt_z\ln z)\right).
$$
Since we must be able to make it arbitrarily negative (by increasing $B$ if necessary) we should require that $t_z\ln z$ is about a constant, say $t_z=1/\ln z$ for $z>e$.
Substituting this into (\ref{tofz}) yields
$$h'\left(1/\ln z\right)=-z\ln z,\quad\mbox{or\quad with $s=1/\ln z$,}\quad h'(s)=-\frac{e^{1/s}}s.$$
Thus we may take
$$h(t)=\int_t^1\frac{e^{1/s}}sds,$$
and we obtain
$$f_M(t)=K\exp\left(-\int_{1-t}^1\frac{e^{1/s}}sds\right),\quad 0<t<1,\quad\mbox{where}\quad K^{-1}=\int_0^1e^{-h(u)}du.$$
(Note that $t_z$ is indeed the local minimum of (\ref{et}).) A graph of the density $f_M$ is given in Figure~\ref{fig:jdens}.

For the lower bound, as
$$p_\delta=K\int_{1-\delta}^1e^{-h(1-t)}dt=K\int^\delta_0e^{-h(t)}dt,$$
we obtain
\begin{eqnarray*}\ln \P(R\ge x)&=&\frac{\ln(1-c)}{\ln(1-cq/x)}\ln\left(K\int_0^{cq/x}e^{-h(t)}dt\right)\\&\sim&-\frac{\ln(1-c)}{cq}x\ln\left(\int_0^{cq/x}e^{-h(t)}dt\right).
\end{eqnarray*}
We need the following lemma which we justify below.
\begin{lem}\lbl{lem:hosp}
\be\lbl{hosp}\frac{y\ln\left(\int_0^{1/y}e^{-h(t)}dt\right)}{e^y}\to-1,\quad\mbox{as}\quad y\to\infty.\ee
\end{lem}
Using this lemma with $y=x/(cq)$ and $c=B/q$ we get, asymptotically,
$$\frac{\ln \P(R\ge x)}{e^{x/B}}\ge-\frac{\ln(1-B/q)}{Be^{x/B}}x\ln\left(\int_0^{B/x}e^{-h(t)}dt\right)\sim
\ln(1-B/q),$$
which implies (\ref{lobdd-j}).\qed

\noindent{\em Proof of Lemma~\ref{lem:hosp}:} We re-write the left--hand side of (\ref{hosp})  as
$$\frac{\ln\left(\int_0^{1/y}e^{-h(t)}dt\right)}{e^y/y},$$
and apply l'Hospital rule. The first differentiation gives
$$\frac{(-1/y^2)e^{-h(1/y)}}{\left(\int_0^{1/y}e^{-h(t)}dt\right)\left(-e^y/y^2+e^y/y\right)}=
\frac{e^{-h(1/y)}e^{-y}/(1-y)}{\int_0^{1/y}e^{-h(t)}dt}.$$
Differentiating again we get
\begin{eqnarray*}&&\frac{(1/y^2)h'(1/y)e^{-h(1/y)}e^{-y}/(1-y)+e^{-h(1/y)}\frac d{dy}\left(e^{-y}/(1-y)\right)}{(-1/y^2)e^{-h(1/y)}}\\ && \quad
=-h'\left(\frac1y\right)\frac{e^{-y}}{1-y}-y^2\frac d{dy}\left(\frac{e^{-y}}{1-y}\right).
\end{eqnarray*}
Since $h'(s)=-e^{1/s}/s$ the first term goes to $-1$ as $y\to\infty$ while the second is $o(1)$.\qed

\vfill
\break

\bibliographystyle{plain}


%

\begin{thebibliography}{10}


\bibitem{air}
G.~Alsmeyer, A.~Iksanov, and U.~R{\"o}sler.
On distributional properties of perpetuities.
{\em J. Theoret. Probab.}, 20:666--682, 2009.

\bibitem{bw}
M.~Bia{\l}kowski and J.~Weso{\l}owski.
\newblock Asymptotic behavior of some random splitting schemes.
\newblock {\em Probab. Math. Statist.}, 22:181--191, 2002.

\bibitem{bond}
L.~Bondesson.
\newblock A general result on infinite divisibility.
\newblock {\em Ann. Probab.}, 7(6):965--979, 1979.

\bibitem{brown}
M.~Brown.
\newblock Error bounds for exponential approximations of geometric
  convolutions.
\newblock {\em Ann. Probab.}, 18(3):1388--1402, 1990.

\bibitem{cl}
J.-F. Chamayou and G.~Letac.
\newblock Explicit stationary distributions for compositions of random
  functions and products of random matrices.
\newblock {\em J. Theoret. Probab.}, 4:3--36, 1991.

\bibitem{eg}
P.~Embrechts and C.~M. Goldie.
\newblock Perpetuities and random equations.
\newblock In {\em Asymptotic statistics (Prague, 1993)}, Contrib. Statist.,
  pages 75--86. Physica, Heidelberg, 1994.

\bibitem{goldie}
C.~M. Goldie.
\newblock Implicit renewal theory and tails of solutions of random equations.
\newblock {\em Ann. Appl. Probab.}, 1(1):126--166, 1991.

\bibitem{gg}
C.~M. Goldie and R.~Gr{\"u}bel.
\newblock Perpetuities with thin tails.
\newblock {\em Adv. in Appl. Probab.}, 28:463--480, 1996.

\bibitem{grey}
D.~R. Grey.
\newblock Regular variation in the tail behaviour of solutions of random
  difference equations.
\newblock {\em Ann. Appl. Probab.}, 4:169--183, 1994.

\bibitem{grin}
A.~K. Grincevi{\v{c}}jus.
\newblock On a limit distribution for a random walk on lines.
\newblock {\em Litovsk. Mat. Sb.}, 15:79--91, 243, 1975.

\bibitem{hm}
P.~Hitczenko and G.~S. Medvedev.
\newblock Bursting oscillations induced by small noise.
\newblock {\em SIAM J. Appl. Math.}, 69:1359 -- 1392, 2009.

\bibitem{jurek}
Z.~J. Jurek.
\newblock Selfdecomposability perpetuity laws and stopping times.
\newblock {\em Probab. Math. Statist.}, 19:413--419, 1999.

\bibitem{kesten}
H.~Kesten.
\newblock Random difference equations and renewal theory for products of random
  matrices.
\newblock {\em Acta Math.}, 131:207--248, 1973.

\bibitem{kn}
M.~Knape and R.~Neininger.
\newblock Approximating perpetuities.
\newblock {\em Methodol. Comput. Appl. Probab.}, 10:507--529, 2008.

\bibitem{kr}
M.~A. Krasnosel{$'$}ski{\u\i} and Ja.~B. Ruticki{\u\i}.
\newblock {\em Convex functions and {O}rlicz spaces}.
\newblock Translated from the first Russian edition by Leo F. Boron. P.
  Noordhoff Ltd., Groningen, 1961.

\bibitem{letac}
G.~Letac.
\newblock A contraction principle for certain {M}arkov chains and its
  applications.
\newblock In {\em Random matrices and their applications (Brunswick, Maine,
  1984)}, number~50 in Contemp. Math., pages 263--273. Amer. Math. Soc.,
  Providence, RI., 1986.

\bibitem{thorinb}
O.~Thorin.
\newblock On the infinite divisibility of the lognormal distribution.
\newblock {\em Scand. Actuar. J.}, (3):121--148, 1977.

\bibitem{thorina}
O.~Thorin.
\newblock On the infinite divisibility of the {P}areto distribution.
\newblock {\em Scand. Actuar. J.}, (1):31--40, 1977.

\bibitem{vervaat}
W.~Vervaat.
\newblock On a stochastic difference equation and a representation of
  nonnegative infinitely divisible random variables.
\newblock {\em Adv. in Appl. Probab.}, 11(4):750--783, 1979.

\bibitem{yan}
N.~Yannaros.
\newblock Randomly observed random walks.
\newblock {\em Comm. Statist. Stochastic Models}, 7(2):219--231, 1991.

\end{thebibliography}

\end{document}